\documentclass[11pt]{article}

\usepackage{amssymb}
\usepackage{epsfig}
\usepackage{hyperref}

\if@twoside \oddsidemargin 5pt \evensidemargin 5pt \marginparwidth 20pt
 \else \oddsidemargin 5pt \evensidemargin 5pt \marginparwidth 20pt \fi

 \mathchardef\hyphen="002D

 \newcommand{\Qed}{\rule{2.5mm}{3mm}}

 \def\NN{\hbox{\sf I\kern-.13em\hbox{N}}}
 \def\RR{\hbox{\sf I\kern-.14em\hbox{R}}}
 \def\ZZ{{\hbox{\sf Z\kern-.43emZ}}}
 \def\CC{\hbox{\sf C\kern -.48emC}}
 \def\QQ{\hbox{\sf C\kern -.48emQ}}
 \def\Cc{\hbox{\sf C\kern -.47em {\raise .48ex \hbox{$\scriptscriptstyle |$}}
   \kern-.5em {\raise .48ex \hbox{$\scriptscriptstyle |$}} }}
 \def\Qq{\hbox{\sf Q\kern -.57em {\raise .48ex \hbox{$\scriptscriptstyle |$}}
   \kern-.55em {\raise .48ex \hbox{$\scriptscriptstyle |$}} }}

\newtheorem{theorem}{Theorem}[section]
\newtheorem{proposition}[theorem]{Proposition}
\newtheorem{lemma}[theorem]{Lemma}
\newtheorem{corollary}[theorem]{Corollary}

\newcommand{\hatr}{half-arc-transitive }
\newcommand{\half}{$\frac{1}{2}$}
\newcommand{\Aut}{\mathrm{Aut}}

\newcommand{\iso}{\cong}

\newcommand{\X}{\mathcal{X}}

\newcommand{\C}{\mathcal{C}}
\newcommand{\G}{\mathcal{G}}
\newcommand{\HH}{\mathcal{H}}
\renewcommand{\P}{\mathcal{P}}
\newcommand{\la}{\langle}
\newcommand{\ra}{\rangle}
\newcommand{\Anc}{\mathrm{Anc}}

\newenvironment{proof}{{\noindent \sc Proof:}}{\hfill $\Qed$}

\begin{document}

\begin{center}
{\bf\large A CLASSIFICATION OF TIGHTLY ATTACHED HALF-ARC-TRANSITIVE GRAPHS OF VALENCY $4$}
\end{center}

\medskip\noindent
\begin{center}
   ~Primo\v z \v Sparl\footnote{{\em E-mail address: Primoz.Sparl@fmf.uni-lj.si}}

\bigskip
{\em  IMFM, University of Ljubljana, Jadranska 19, 1111 Ljubljana, Slovenia}
\end{center}




\vskip 1cm

\hrule
\begin{abstract}
A graph
is said to be {\em half-arc-transitive} if its automorphism group acts
transitively on the set of its vertices and edges but not on the set of its arcs.
With each half-arc-transitive graph of valency $4$ a collection of the so called {\em alternating cycles}
is associated, all of which have the same even length. Half of this length is called the {\em radius}
of the graph in question.
Moreover, any two adjacent alternating cycles 
 have the same number of common vertices. If this number, the so called 
{\em attachment number}, coincides with the radius,
we say that the graph is {\em tightly attached}.
In {\em J.~Combin.~Theory Ser.~B} {73} (1998) 41--76, Maru\v si\v c gave a classification of
tightly attached \hatr graphs of valency $4$ with odd radius.
In this paper the even radius tightly attached graphs of valency $4$ are classified, thus 
completing the classification of all tightly attached half-arc-transitive graphs
of valency $4$.
\end{abstract}

\begin{quotation}
\end{quotation}
\hrule\medskip


\section{Introductory remarks}
\indent

Throughout this paper graphs are assumed to be finite and,
unless stated otherwise, simple, connected and undirected.
For group-theoretic concepts not defined here we
refer the reader to \cite{BW, DM, W}, and for graph-theoretic terms not 
defined here we refer the reader to \cite{BM}. In this paper we let $\ZZ_n$ denote the ring of
residue classes modulo $n$ and we let $\ZZ_n^*$ denote the set of invertible elements of $\ZZ_n$. 
At times it will be convenient to view elements of $\ZZ_n$ as integers,
for instance if $\rho$ is an element of some group with $\rho^n = 1$ and if $r \in \ZZ_n$, we
let $\rho^r$ represent $\rho^k$ for any $k$ in the equivalence class $r$. This should cause no confusion.
For basic notation and other conventions see Section~\ref{sec:thegraphs}.

Let $X$ be a graph. We let $V(X)$, $E(X)$ and $A(X)$ denote the set of 
vertices, edges and arcs of $X$, respectively.
The graph $X$ is said to be {\em vertex-transitive}, {\em edge-transitive} and 
{\em arc-transitive} provided that its automorphism group $\Aut X$ acts transitively on the
set of its vertices, edges and arcs, respectively. Moreover,  $X$ is said to be
{\em half-arc-transitive} if it is vertex- and edge- but not arc-transitive.
More generally, by a {\em half-arc-transitive} action of a 
subgroup $G \leq \Aut X$ on $X$ we mean a vertex- and edge- but not
arc-transitive action of $G$ on $X$.  
In this case we say that  $X$ is  $G$-{\em half-arc-transitive}.
As demonstrated in \cite[7.53, p.~59]{T66} by Tutte, the valency of a graph admitting 
a half-arc-transitive group action is necessarily even.
A few years later Tutte's question as to the existence  of half-arc-transitive graphs
of a given even valency was answered by Bouwer~\cite{B70} 
with a construction of a $2k$-valent half-arc-transitive graph for every $k\geq2$.
The smallest graph in Bouwer's family has 54 vertices and
valency 4. Doyle~\cite{PD76} and Holt~\cite{H81} independently found one with 27 vertices,
a graph that is now known to be the smallest half-arc-transitive graph~\cite{AMN94}.

Interest in the study of these graphs
reemerged in the nineties with a series of papers
dealing mainly with classification of certain restricted classes of such graphs 
as well as with various methods of constructions of new families of such graphs
\cite{AMN94,AX94,MX97,TX94,TW89,Wa94,Xu92}.
These graphs have remained an active topic of research to this day,
with a number of papers dealing with their structural properties;
see \cite{CM03,DX98,LiS01,MM99,MM02,DM98,MNmas,MNgrt,MNjams}.
However, graphs admitting  half-arc-transitive group actions
are in a 1-1-correspondence with the so called orbital graphs of groups with non-self-paired orbitals.
(Given a transitive permutation group $G$
acting on a set $V$, let $\mathcal{O}$ be a nontrivial, that is $\mathcal{O} \neq \{(v,v)\ |\ v \in V\}$,
and non-self-paired orbital, that is different from its paired orbital $\mathcal{O}^T = \{(u,v)\ |\ (v,u) \in \mathcal{O}\}$,
in the natural action of $G$ on $V \times V$. 
Then the graph with vertex set  $V$ and edge set $\{ uv\ |\ (u,v) \in \mathcal{O}\}$ is half-arc-transitive.
Conversely, every graph admitting a half-arc-transitive group action arises in this way.)
The  classification of the whole class of  half-arc-transitive graphs is therefore presently beyond our reach, and it thus
seems only natural to restrict our consideration to some special classes of these  graphs.
There are several approaches that are currently being taken, such as for example,
investigation of (im)primitivity of half-arc-transitive group actions on graphs 
\cite{DX98, LLM04, TX94},
geometry related questions about half-arc-transitive graphs  \cite{DAN04, MN98, MPi99},
and questions concerning  classification for various restricted classes of half-arc-transitive graphs
\cite{AX94, DM98, Wa94, Xu92},  to mention just a few.

In view of the fact that $4$ is the smallest admissible valency
for a half-arc-transitive graph, 
special attention has rightly been given to 
the study of half-arc-transitive graphs of valency $4$.
However, even this restricted class of graphs is very rich
and only partial results have been obtained thus far.
One of the possible approaches is the study of these graphs 
via the corresponding vertex stabilizers with some promising  results proved in \cite{MNgrt}.
An alternative point of view, more geometric in nature,
was first presented in \cite{DM98}.
The idea is to obtain some insight into structural properties of the graph by studying its so called alternating cycles.
We give a brief explanation of the concepts involved below.

Let $X$ be a $G$-half-arc-transitive graph of valency $4$ where $G \leq \Aut X$ and let $D(X)$ be one
of the two oriented graphs corresponding to this half-arc-transitive action of $G$,
obtained by orienting an arbitrary edge in one of the two possible ways and then 
applying the action of $G$ to obtain a unique orientation of the whole edge set of $X$.
We say that a cycle $C$ in $X$ of even length is a $G$-{\em alternating cycle} if
its vertices are alternately the heads and the tails (in $D(X)$) of their two incident edges in $C$.
It was proved in \cite[Proposition 2.4.]{DM98} that all $G$-alternating cycles of $X$ have
equal length $2r_G(X)$ for some $r_G(X) \geq 2$. The parameter $r_G(X)$ is called the $G$-{\em radius} of $X$. 
We say that the $G$-alternating cycles $C_1$ and $C_2$ are adjacent if they share a common vertex. Since $G$
acts transitively on the vertices of $X$, all pairs of adjacent $G$-alternating cycles in $X$ have an 
equal number of common vertices.
We call this number the $G$-{\em attachment number} of $X$ an denote it by $a_G(X)$.
Note that the sets of common vertices of $G$-alternating cycles, called the $G$-{\em attachment sets}, are
blocks of imprimitivity for $G$.
The relation of the parameters $r_G(X)$ and $a_G(X)$ is important. 
It was shown in \cite{DM98} and \cite{MP99}
that $a_G(X)$ divides $2r_G(X)$ and is at most $r_G(X)$ in the case when $X$ is \hatr. 
If $a_G(X) = r_G(X)$, we say that
$X$ is $G$-{\em tightly attached}. At the other extreme, we say that $X$ is 
$G$-{\em loosely attached} if $a_G(X) = 1$, and we say that $X$ is $G$-{\em antipodaly attached} if
$a_G(X) = 2$. The importance of these three families of graphs is suggested by results from 
\cite{MP99}, where among other, it was shown that  every $G$-half-arc-transitive graph of valency $4$ is either $G$-tightly attached
or it is a cover either of a $G$-loosely attached or of a $G$-antipodaly attached graph. In all of the above 
terminology the prefix $G$ is omitted when $G = \Aut X$.
Let us also mention that infinite families of half-arc-transitive graphs with prescribed attachment numbers 
were constructed in \cite{MW00}.

As the structure of tightly attached graphs
seems most natural and easy to understand, the first step in the classification of the half-arc-transitive graphs
of valency $4$ is thus to classify these graphs. In 1998 Maru\v si\v c gave a complete classification
of the odd radius graphs. His result is the following

\begin{theorem}[{\cite[Theorem~3.4.]{DM98}}]
\label{the:oddtheorem}
A graph $X$ is a tightly attached half-arc-transitive graph of valency $4$ and
odd radius $n$ if and only if $X \iso \X_o(m,n;r)$, where $m \geq 3$ and 
$r \in \ZZ_n^*$ satisfies $r^m =\pm 1$, and moreover none of the following
conditions is fulfilled:
\begin{itemize}
	\item[(i)] $r^2 = \pm 1$;
	\item[(ii)] $(m,n;r) = (3,7;2)$;
	\item[(iii)] $(m,n;r) = (6,7n_1;r)$, where $n_1 \geq 1$ is odd and coprime to $7$, $r^6 = 1$, and
	there exists a unique solution $r' \in \{r, -r, r^{-1}, -r^{-1}\}$ of the equation
	$2-r-r^2=0$ such that $7(r'-1)=0$ and $r' \equiv 5 \pmod 7$.
\end{itemize}
\end{theorem}
For the definition of the graphs $\X_o(m,n;r)$ see Section~\ref{sec:thegraphs}.

In \cite{MP99} the graphs of valency $4$ which
admit a half-arc-transitive group action relative to which the graph is tightly attached and 
has even radius were classified by Maru\v si\v c and Praeger. But the question of which of these graphs are indeed 
half-arc-transitive and which are arc-transitive was not answered. In $2004$ Wilson~\cite{SW04}
found an alternative way of describing the graphs of valency $4$ which
admit a half-arc-transitive group action relative to which the graph is tightly attached.
He showed that these graphs are the so called power spider and mutant power spider graphs. But even with this
improvement the question of half-arc-transitivity of these graphs remained unsolved. 

It is the aim of this paper to
resolve this question. We improve the  results of Maru\v si\v c and Praeger on
graphs of valency $4$ admitting a half-arc-transitive group action relative to which the graph is tightly
attached of even radius and then determine
precisely which of these graphs are \hatr and which are arc-transitive. 
Together with the above mentioned classification of the odd radius case, this gives the 
complete classification of tightly-attached \hatr graphs of valency $4$. Our main result is
the following

\begin{theorem}
\label{the:thetheorem}
A graph $X$ is a tightly attached half-arc-transitive graph of valency $4$ and
even radius $n$ if and only if $X \iso \X_e(m,n;r,t)$, where $m \geq 4$ is even,
$r \in \ZZ_n^*$, $t \in \ZZ_n$ are such that 
$r^m = 1$, $t(r-1) = 0$ and $1 + r + \cdots + r^{m-1} + 2t = 0$, and none of the following
two conditions is fulfilled:
\begin{itemize}
	\item[(i)] $r^2 = \pm 1$;
	\item[(ii)] $m = 6$, $n = 14n_1$, where $n_1$ is coprime to $7$, precisely
	one of $\{r, -r, r^{-1}, -r^{-1}\}$ solves the equation $2-x-x^2=0$ and if we let $r'$ be this unique
	solution and let $t' = t$ in case $r' \in \{r,r^{-1}\}$ and $t' = t + r + r^3 + \cdots + r^{m-1}$ in case
	$r' \in \{-r, -r^{-1}\}$, then $r' \equiv 5 \pmod 7$,
	$7(r'-1) = 0$, $t' \equiv 0 \pmod 7$ and $2+r'+t'=0$.
\end{itemize}
\end{theorem}
For the definition of the graphs $\X_e(m,n;r,t)$ see Section~\ref{sec:thegraphs}. 
In Proposition~\ref{pro:alliso}
we determine precisely which pairs of the half-arc-transitive graphs $\X_e(m,n;r,t)$ are isomorphic.
Along the way we also complete the work of \v Sajna (see \cite{MS98}) in determining which of the 
metacirculants $M(r;4,n)$ are half-arc-transitive. For details see Section~\ref{sec:proof}.

The paper is organized as follows. In Section~\ref{sec:thegraphs} we introduce the graphs $\X_e(m,n;r,t)$.
We then prove that a graph of valency $4$ admitting a 
half-arc-transitive group action relative to which the graph is tightly attached of even radius is
isomorphic to some $\X_e(m,n;r,t)$. The investigation of whether or not such a graph is half-arc-transitive
is based on the ideas introduced in \cite{DM98} where the odd radius graphs were classified. The idea is to
determine the possible $8$-cycles of the graph and then investigate the
interplay of these $8$-cycles with the $2$-paths of the graph. It is this information that gives an
insight into the arc- or half-arc-transitivity of the graphs in question. The terminology
and basic properties are introduced in Section~\ref{sec:notation} and the possible $8$-cycles are
investigated in Section~\ref{sec:8cycles}. In the subsequent sections we then
investigate the above mentioned interplay of $8$-cycles and $2$-paths 
depending on the number of alternating cycles the graph in
question has. The proof of our main result, Theorem~\ref{the:thetheorem}, is then laid out in 
Section~\ref{sec:proof}. Proposition~\ref{pro:alliso} and the results concerning metacirculants
$M(r;4,n)$ are also stated and proved there.


\section{The even radius graphs}
\label{sec:thegraphs}

\indent
Let $X$ be a $G$-half-arc-transitive graph for some $G \leq \Aut X$ and let $D(X)$ be one of the two
oriented graphs corresponding to this action of $G$.
Let $u,v \in V(X)$ be adjacent (we denote this by $u \sim v$ and we
denote the corresponding edge by $uv$). Then of course
either $(u,v)$ or $(v,u)$ is an arc of $D(X)$. In the former case we say that $u$ is the {\em tail}
and $v$ is the {\em head} of $(u,v)$, and we say
that $u$ is the {\em predecessor} of $v$ and $v$ is the {\em successor} of $u$. 
Throughout the paper we shall constantly be switching from the two viewpoints regarding $X$, 
namely as an undirected graph or an oriented graph $D(X)$. This should cause no confusion.

For the sake of completeness we define the graphs from Theorem~\ref{the:oddtheorem}.
For each $m \geq 3$, $n \geq 3$ odd, $r \in \ZZ_n^*$, where $r^m = \pm 1$, and $t \in \ZZ_n$
let $\X_o(m,n;r)$ be the graph with vertex 
set $V = \{u_i^j\ |\ i \in \ZZ_m,\ j \in \ZZ_n\}$ 
and edges defined by the following adjacencies:
$$ u_i^j \sim u_{i+1}^{j \pm r^i}\quad ;\quad i \in \ZZ_n,\ j \in \ZZ_m.$$
(The subscript $o$ in the symbol $\X_o(m,n;r)$ is meant to 
indicate that $n$ is an odd integer.) Note that the 
graphs $\X_o(m,n;r)$ correspond to the graphs $X(r;m,n)$ introduced in \cite{DM98}.

We now review and improve somewhat the results of Maru\v si\v c and Praeger from \cite{MP99}.
Let $X$ be a graph of valency $4$ admitting a half-arc-transitive subgroup $G$ of $\Aut X$.
Suppose $X$ is tightly $G$-attached with an even $G$-radius $n \geq 4$.
Note that we need not consider the case $n = 2$, as  $X$ is a lexicographic product of a cycle by $2K_1$
and thus arc-transitive in that case.
Let $m$ denote the number of $G$-alternating cycles of $X$ and let $\Sigma$ denote
the set of $G$-attachment sets of $X$ (recall that $\Sigma$ is a complete imprimitivity block system for $G$).
Let now $C = v_0 v_1 \ldots  v_{2n-1}$ be
any $G$-alternating cycle of $X$. 
It may be seen that  there exists some $\rho \in G$ whose restriction to $C$ is 
$(v_0v_2\ldots v_{2n-2})(v_1v_3\ldots v_{2n-1})$  (see \cite{MP99}).
Since $n \geq 4$
and $X$ is tightly $G$-attached, \cite[Lemma~3.5]{MP99} implies that $G_{v_0} \iso \ZZ_2$. 
Let $\tau \in G$ be the unique nonidentity element of $G_{v_0}$. 
Furthermore, let $\sigma \in G$ be such that 
$v_0 \sigma = v_{2n-1}$. (Note that in \cite{MP99} $\sigma$ was chosen so as to map $v_0$ to $v_1$
but we prefer this choice in order to obtain a more convenient description of $X$.) Clearly, $\sigma$ cyclically permutes the
$m$ blocks of $\Sigma$. Let $K$ denote the kernel of the action of $G$ on $\Sigma$. 
Then the following theorem, which is a slight improvement of \cite[Theorem~4.2]{MP99}, holds.

\begin{theorem}
\label{the:evengroup}
The permutations  $\rho, \sigma, \tau $  generate $G$ and satisfy the following relations
\begin{equation}
\label{eq:relations}
\rho^n = \tau^2 = 1,\ \sigma^m = \rho^t,\ \rho^\tau = \rho^{-1},\ \rho^\sigma = \rho^r,\ \tau^\sigma = \tau\rho^{-1}
\end{equation}
where both $n \geq 4$ and $m \geq 3$ are even
and $r \in \ZZ_n^*$, $t \in \ZZ_n$ are such that
\begin{equation}
\label{eq:rt}
r^m = 1,\quad t(r-1) = 0\quad \mathrm{and}\quad 1+r+r^2+\cdots +r^{m-1} + 2t = 0.
\end{equation}
Further, $K = \la \rho, \tau \ra \iso D_{2n}$.
\end{theorem}

\begin{proof}
Following the proof of \cite[Theorem~4.2]{MP99} we find that $K = \la \rho, \tau \ra = D_{2n}$,
$G = \la \rho, \sigma, \tau \ra$ has order $2mn$ and $\rho$ is of order $n$. 
Moreover, $\sigma^m = \rho^t$ for some $t \in \ZZ_n$ and 
there exists some $r \in \ZZ_n^*$ such that $\rho^{\sigma} = \rho^r$,
$r^m = 1$ and $t(r-1) = 0$. Furthermore,   there
exists some $k \in \ZZ_n$ such that $\tau^\sigma = \tau\rho^k$. 

Since $\tau \in G_{v_0}$ is nontrivial, it interchanges $v_1$ and $v_{2n-1}$. Thus
$v_{2n-1}\tau^\sigma = v_0\sigma\tau^{\sigma} = v_0\tau\sigma = v_0\sigma$. 
On the other hand,
$v_{2n-1}\tau^{\sigma} = v_{2n-1}\tau\rho^k = v_1\rho^k = v_0\sigma \rho \rho^k = v_0\sigma \rho^{k+1}$.
Therefore $v_0\sigma = v_0\sigma\rho^{k+1}$, and so $\rho^{k+1} \in G_{v_{2n-1}}$.
Then clearly $\rho^{k+1} = 1$, implying $k = n-1$. Therefore, $\tau^{\sigma} = \tau\rho^{-1}$, and so equations
$\tau^{\sigma^m} = \tau\rho^{-(1+r+\cdots + r^{m-1})}$ and $\tau^{\sigma^m} = \tau^{\rho^t} = \tau\rho^{2t}$
give us $1+r+\cdots + r^{m-1} + 2t = 0$. Combining this with  the fact that $n$ is even 
we get that $1+r+\cdots + r^{m-1} \equiv 0 \pmod 2$. But since $r$ is coprime to $n$ and hence is odd, 
it follows that $m$ is even. 
\end{proof}\bigskip

We now introduce a family of graphs that will play a central role  in this paper.
For all even integers $m \geq 4$, $n \geq 4$ and for each $r \in \ZZ_n^*$, $t \in \ZZ_n$ satisfying
\begin{equation}
\label{eq:Xparamcond}
	r^m = 1, \quad t(r-1) = 0\quad \mathrm{and}\quad 1 + r + \cdots + r^{m-1} + 2t = 0,
\end{equation}
let $\X_e(m,n;r,t)$ be the graph with vertex 
set $V = \{u_i^j\ |\ i \in \ZZ_m,\ j \in \ZZ_n\}$ 
and edges defined by the following adjacencies:
$$ u_i^j \sim \left\{\begin{array}{lll}
	u_{i+1}^j,\ u_{i+1}^{j + r^i} & ; & i \in \ZZ_m \setminus \{m-1\},\ j \in \ZZ_n \\ \\
	u_{0}^{j+t},\ u_0^{j+r^{m-1}+t} & ; & i = m-1,\ j \in \ZZ_n .\end{array}\right. $$
(The subscript $e$ in the symbol $\X_e(m,n;r,t)$ is meant to 
indicate that $n$ is an even integer.)
Let $\rho$, $\sigma$ and $\tau$ be the permutations defined on $V$
by the following rules
\begin{equation}
\label{eq:rho}
	u_i^j\rho = u_i^{j+1} \quad  ; \quad i \in \ZZ_m,\  j \in \ZZ_n 
\end{equation}
\begin{equation}
\label{eq:sigma}
	u_i^j\sigma = \left\{\begin{array}{lll}
	u_{i+1}^{rj} & ; & i \in \ZZ_m \setminus \{m-1\},\ j \in \ZZ_n\\ \\
	u_{0}^{rj+t} & ; & i = m-1,\ j \in \ZZ_n\end{array}\right.
\end{equation}
\begin{equation}
\label{eq:tau}
	u_i^j\tau = \left\{\begin{array}{lll}
	u_{0}^{-j} & ; & i = 0,\ j \in \ZZ_n\\ \\
	u_{i}^{1+r+\cdots +r^{i-1} - j} & ; & i \in \ZZ_m\setminus \{0\},\ j \in \ZZ_n .\end{array}\right. 
\end{equation}
Clearly $\rho$ and $\sigma$ are automorphisms of $\X_e(m,n;r,t)$ since $t(r-1)=0$.
As for $\tau$, it is clear that every edge connecting vertices with
subscripts $i$ and $i+1$ is mapped to an edge when $i \neq m-1$. Moreover,
the neighbors $u_0^{j+t}$ and $u_0^{j+r^{m-1}+t}$ of $u_{m-1}^j$ are mapped to
$u_0^{-j-t}$ and $u_0^{-j-r^{m-1}-t}$, respectively. Thus since, in view of $(\ref{eq:Xparamcond})$, we have 
$u_{m-1}^{j}\tau = u_{m-1}^{1+r+\cdots + r^{m-2} - j} = u_{m-1}^{-j-r^{m-1}-2t}$, $\tau$
is also an automorphism of $\X_e(m,n;r,t)$. This implies that $H = \langle \rho, \sigma, \tau\rangle$ acts
half-arc-transitively on   $\X_e(m,n;r,t)$.

Using Theorem~\ref{the:evengroup} together with the proof of \cite[Theorem~4.5]{MP99}  
it may be seen that the converse also holds.

\begin{theorem}
\label{the:Xgraphs}
Let $X$ be a graph of valency $4$ admitting a half-arc-transitive subgroup $G$ of $\Aut X$ 
relative to which it is tightly attached of even radius.
Then $X$ is isomorphic to $\X_e(m,n;r,t)$ for some even $m \geq 4$,
and some $r \in \ZZ_n^*$ and $t \in \ZZ_n$ satisfying (\ref{eq:Xparamcond}).
\end{theorem}

This theorem thus classifies all graphs of valency $4$ which admit a \hatr subgroup of automorphisms
relative to which the graph is tightly attached with even radius. In the rest of this paper we determine
which of the graphs $\X_e(m,n;r,t)$ are \hatr and which are arc-transitive. A complete
classification of the tightly attached \hatr graphs of even radius and valency $4$ is thus obtained.


\section{Notation and preliminary results}
\label{sec:notation}

\indent
In this section we let $m,n \geq 4$ be even and we let $r\in \ZZ_n^*$, $t\in \ZZ_n$
satisfy (\ref{eq:Xparamcond}). 
We use $X$ as a shorthand notation for the graph $\X_e(m,n;r,t)$ and we let
$\rho$, $\sigma$ and $\tau$ be as in (\ref{eq:rho}), (\ref{eq:sigma}) and (\ref{eq:tau}),
respectively.
We let $H = \langle \rho, \sigma , \tau \rangle$ and we let 
$X_i = \{u_i^j\ |\ j \in \ZZ_n\}$, $i \in \ZZ_m$, denote the orbits of $\rho$.
Clearly, the sets $X_i$, $i \in \ZZ_m$, are the attachment sets in the
half-arc-transitive action of $H$ on $X$, and, of course, blocks of imprimitivity for $H$.

\begin{proposition}
\label{pro:r^2}
Let $m,n \geq 4$ be even integers and let $r\in \ZZ_n^*$, $t \in \ZZ_n$ satisfy
(\ref{eq:Xparamcond}). If $r^2 = \pm 1$ then $\X_e(m,n;r,t)$ is arc-transitive.
\end{proposition}

\begin{proof}
Let $X = \X_e(m,n;r,t)$. As noted above,
$H$ acts half-arc-transitively on $X$. We thus only need to show that
there exists an automorphism of $X$ interchanging two adjacent vertices of $X$. 
We distinguish two cases depending on whether $r^2 = 1$ or $r^2 = -1$.

Suppose first that $r^2 = 1$. This implies that $r^i = r$ for $i$ odd and $r^i = 1$ for $i$ even.
In particular $r^{m-1} = r$ since $m$ is even.
Let $\varphi$ be the permutation of $V(X)$ defined by the 
rule
$$ u_i^j\varphi = \left\{\begin{array}{lll} u_0^{-rj} & ; & i = 0,\ j \in \ZZ_n \\ 
			u_{m-i}^{-rj-t} & ; & i \in \ZZ_m\setminus \{0\},\ j \in \ZZ_n.\end{array}\right.
$$
Note that since $r \in \ZZ_n^*$, $\varphi$ is indeed a permutation of $V(X)$.
We claim that $\varphi$ is in fact an automorphism of $X$. To see this we show that
each edge joining a vertex in  $X_i$ with a vertex in  $X_{i+1}$ is mapped to an edge. 
For instance, if $i = 0$ then 
for any $j \in \ZZ_n$ we have $u_{0}^j \sim u_1^j, u_1^{j+1}$ and the images of 
these three vertices under $\varphi$ are 
$u_0^{-rj}$, $u_{m-1}^{-rj-t}$ and $u_{m-1}^{-rj-r-t}$, respectively, and so
$u_0^j\varphi \sim u_1^j\varphi, u_1^{j+1}\varphi$. 
The arguments for $0 < i < m-1$ only depend on parity of
$i$ and are left to the reader. Finally, for any $j \in \ZZ_n$ we have
$u_{m-1}^j \sim u_0^{j+t}, u_0^{j+r+t}$ and the images of these three vertices under $\varphi$
are $u_{1}^{-rj-t}$, $u_0^{-rj-t}$ and $u_{0}^{-rj-1-t}$, respectively (recall that $rt = t$), so
$\varphi$ is indeed an
automorphism of $X$. Thus $X$ is arc-transitive since $\varphi\sigma$ interchanges 
adjacent vertices $u_0^0$ and $u_1^0$.

Suppose now that $r^2 = -1$. Then
$$ r^i = \left\{\begin{array}{lll} 	1 & ; & i \equiv 0 \pmod 4\\
					r & ; & i \equiv 1 \pmod 4\\
					-1& ; & i \equiv 2 \pmod 4\\
					-r& ; & i \equiv 3 \pmod 4.\end{array}\right.$$
Note that $m \equiv 0\pmod 4$ in view of $r^m = 1$, and so $r^{m-1} = -r$.
Let $\psi$ be the permutation of $V(X)$ defined by the rule
$$ u_i^j\psi = \left\{\begin{array}{lll}
	u_{0}^{-rj} 	& ; 	& i = 0,\ j \in \ZZ_n\\   	
	u_{m-i}^{-rj-t} & ; 	& i \in \ZZ_m \setminus \{0\},
	\ i \equiv 0 \pmod 4\ \mathrm{or}\ i \equiv 3 \pmod 4,\ j \in \ZZ_n \\
	u_{m-i}^{-rj+r-t} & ; 	& i \in \ZZ_m,\ 
	i \equiv 1 \pmod 4\ \mathrm{or}\ i \equiv 2 \pmod 4,\ j \in \ZZ_n.\end{array}\right.$$
Using the fact that $r^{m-i} = r^i$ if $i \equiv 0 \pmod 2$ and
$r^{m-i} = -r^i$ if $i \equiv 1 \pmod 2$,
it is easy to see that $\psi \in  \Aut X$.
We leave the details to the reader. Since $\psi\tau\sigma$ interchanges adjacent vertices $u_0^0$ and
$u_1^0$, the graph $X$ is arc-transitive.
\end{proof}\bigskip
    
Our approach in determining whether $X$ is \hatr or arc-transitive when $r^2 \neq \pm1$
relies on a thorough analysis of $8$-cycles in $X$ and their interplay with $2$-paths of $X$,
an idea used also in the classification of odd radius half-arc-transitive graphs of valency $4$
in \cite{DM98}. The terminology too,  basically follows that of \cite{DM98}.

\begin{figure}
\begin{center}
\includegraphics[scale=0.5]{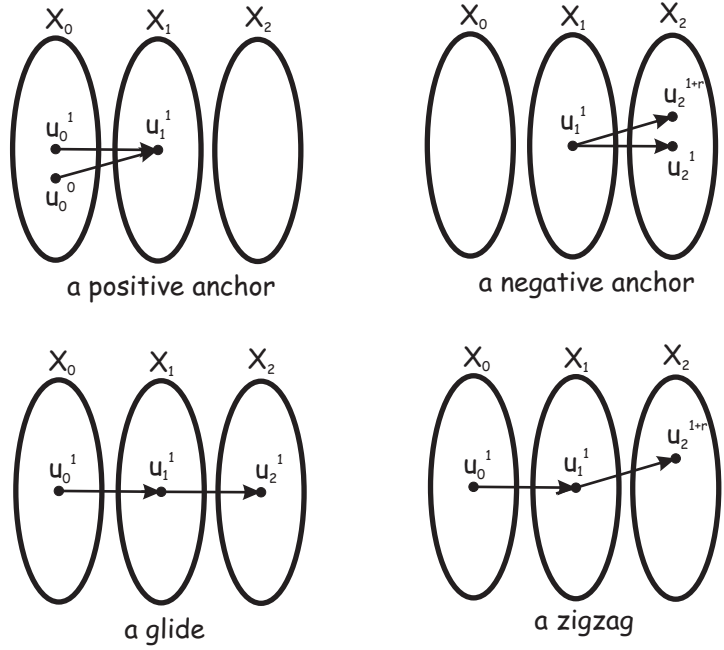}
\caption{The representatives of $H$-orbits of $2$-paths.}
\label{fig:2-paths}
\end{center}
\end{figure}

Note that the group $H$ has four orbits
in its natural action on the set of $2$-paths of $X$. Following \cite{DM98} we introduce the
notation concerning $2$-paths
of different $H$-orbits as follows. We call any $2$-path
whose endvertices belong to the same set $X_i$ an {\em anchor}. Note that the group $H$
has two orbits on the set $Anc X$ of all anchors of $X$, namely one
containing the anchor $u_0^1u_1^1u_0^0$ and one containing the anchor $u_2^{1+r}u_1^1u_2^1$.
We denote these two orbits by $Anc^+X$ and $Anc^-X$, respectively (see Figure~\ref{fig:2-paths}).
The anchors in $Anc^+X$ will be called  {\em positive anchors} and the
anchors in $Anc^-X$ will be called  {\em negative anchors}. The group $H$ has two additional
orbits on the set of all $2$-paths of $X$. The first one  contains the $2$-path $u_0^1u_1^1u_2^1$,
is denoted by $Gli X$,  and its elements are referred to as {\em glides}.
The second one  contains the $2$-path  $u_0^1u_1^1u_2^{1+r}$, is denoted
by $Zig X$, and its elements are referred to as  {\em zigzags}. 
Note that there is precisely one positive and precisely
one negative anchor having a given vertex as its internal vertex. On the other hand there are precisely two
glides and precisely two zigzags having a given vertex as its internal vertex. Thus we have 
\begin{equation}
\label{eq:orblen}
	|Anc^+X| = |Anc^-X| = mn\quad \mathrm{and}\quad |Gli X| = |Zig X| = 2mn.
\end{equation}
Note also that $u_{m-1}^ju_0^{j+t}u_1^{j+t}$ is a glide and $u_{m-1}^ju_0^{j+t}u_1^{j+1+t}$ is a zigzag.

The next proposition links the problem of determining whether or not $X$ is half-arc-transitive to
the investigation of the action of its automorphism group on the set of $2$-paths of $X$.

\begin{proposition}
\label{pro:arc-tr}
Let $m,n \geq 4$ be even integers and let $r\in \ZZ_n^*$, $t \in \ZZ_n$ be such that 
(\ref{eq:Xparamcond}) holds.  Then $X = \X_e(m,n;r,t)$ is arc-transitive if and only if
either $r^2 = \pm 1$ or the automorphism group $\Aut X$ does not fix the set of anchors $Anc X$.
\end{proposition}

\begin{proof}
Clearly if some $\varphi \in \Aut X$ does not fix the set of anchors, then for some $\psi \in H$ the
automorphism $\varphi \psi$ interchanges a pair of
adjacent vertices of $X$, and so $X$ is arc-transitive. In view of Proposition~\ref{pro:r^2} we thus only
need to show that if $X$ is arc-transitive and $\Aut X$ fixes the set $\Anc X$, then $r^2 = \pm 1$.
Now if $X$ is arc-transitive then there exists some $\varphi \in \Aut X$ fixing $u_0^0$ and
mapping $u_1^0$ to $u_{m-1}^{-t}$. Since $\varphi$ maps the anchor $u_1^1u_0^0u_1^0$ to an anchor, 
we have
$u_1^1\varphi = u_{m-1}^{-r^{m-1}-t}$. Similarly the anchor $u_0^0u_1^1u_0^1$ is mapped to an anchor, so
$u_0^1\varphi = u_0^{-r^{m-1}}$. Continuing this way we get that 
$u_1^r\varphi = u_{m-1}^{-t-rr^{m-1}} = u_{m-1}^{-t-1}$. Note further that since $u_2^0u_1^0u_2^{r}$
is an anchor and $u_0^0\varphi = u_0^0$, 
we have $\{u_2^0, u_2^{r}\}\varphi = \{u_{m-2}^{-t}, u_{m-2}^{-t-r^{m-2}}\}$. 

But $u_1^r \sim u_2^r$, and so $u_{m-1}^{-t-1}$ is adjacent to one of the vertices $u_{m-2}^{-t}, u_{m-2}^{-t-r^{m-2}}$.
If $u_{m-1}^{-t-1} \sim u_{m-2}^{-t}$, then $r^{m-2} = -1$, and so $r^2 = -1$. 
If on the other hand, $u_{m-1}^{-t-1} \sim u_{m-2}^{-t-r^{m-2}}$, then $-1 = -r^{m-2}$, and  so $r^2 = 1$,
completing the proof of Proposition~\ref{pro:arc-tr}.
\end{proof}\bigskip

Continuing with further notation and terminology, we let $W$ be a simple walk of length
$d$ in $X$. To each internal vertex $v$ of $W$ we assign one of the symbols $a,g$ or $z$, depending
on whether the corresponding $2$-path of $W$ having $v$ as its internal vertex 
is an anchor, a glide or a zigzag, respectively. In this way
a sequence of symbols from the set $\{a,g,z\}$ is assigned to $W$. If $W$ is a cycle, then every vertex of $W$ is
internal so that the length of the obtained sequence is $d$. Otherwise the length of the sequence is $d-1$.
We let the equivalence class of all sequences obtained from the above sequence by a reflection or a cyclic rotation
in case when $W$ is a cycle, and just by a reflection in case when $W$ is a path, be the {\em code} of $W$. We let
the {\em refinement of the code} of $W$ be the sequence obtained from the code $C$ of $W$ by replacing each $a$ in $C$ by
$a^+$ or $a^-$ depending on whether $a$ corresponds to a positive or a negative anchor. The {\em trace} of $W$ is obtained
from its code by replacing each $g$ and $z$ by an $n$. Therefore the trace of a walk distinguishes solely between
anchors and nonanchors.  

We say that a
cycle of length $d$ of $X$ is {\em coiled} if its trace is $n^d$ and is {\em noncoiled} otherwise. 
The next observation, essentially  a translation of \cite[Proposition 4.5]{DM98} for
even radius graphs, is self-explanatory.

\begin{proposition}
\label{pro:anchors}
Let $C$ be a noncoiled cycle of $\X_e(m,n;r,t)$, where $m,n \geq 4$ are even integers and $r\in \ZZ_n^*$, 
$t \in \ZZ_n$ are such that (\ref{eq:Xparamcond}) holds. Then 
positive and negative anchors alternate on $C$.
\end{proposition}

\begin{proposition}
\label{pro:nonanchors}
Let $C$ be a cycle of $\X_e(m,n;r,t)$, where $m,n \geq 4$ are even integers and $r\in \ZZ_n^*$, 
$t \in \ZZ_n$ are such that (\ref{eq:Xparamcond}) holds. Then $C$ has an
even number of glides and an even number of zigzags.
\end{proposition}

\begin{proof}
Note that since $m$ is even, $C$ is of even length.
Recall that $D(X)$ is one of the two oriented graphs corresponding to
the half-arc-transitive action of $H$ on $X$. 
Observe that the group $\langle \rho, \sigma \rangle$ has precisely two orbits
on the set of arcs of $D(X)$, one orbit, denoted by $\mathcal{O}_1$, corresponding to the edge $u_0^0u_1^0$ and the
other, denoted by $\mathcal{O}_2$, corresponding to the edge $u_0^0u_1^1$. 
We now assign a binary sequence to $C$ by assigning  a symbol $0$ or $1$
to each edge of $C$ depending on whether it is in $\mathcal{O}_1$ or $\mathcal{O}_2$, respectively. 

Note that if two consecutive symbols in this sequence are equal, then
the common vertex of the two arcs in question is the internal vertex of a glide,
and is an internal vertex of anchor or a zigzag otherwise. 
By Proposition~\ref{pro:anchors}, $C$ has an even number of anchors, 
and furthermore, in a cyclic traversing of the above mentioned binary sequence, 
the number of times we get different consecutive symbols is even.
Consequently, the number of zigzags of $C$ is even, and since
the length of $C$ is even, the number of glides is also even.
\end{proof}\bigskip

Let $\C$ be a union of $H$-orbits of cycles of $X$ and let $P$ be a path of $X$. The number of cycles of 
$\C$ containing $P$ as a subgraph will be called the {\em $\C$-frequency} of $P$ and will be denoted by $v(\C,P)$.
In particular if $d \geq 3$ and $\C$ is the set of all $d$-cycles of $X$ then we let the {\em $d$-frequency}
$v(d,P)$ of $P$ be $v(\C,P)$.

\begin{proposition}
\label{pro:equalfreq}
Let $X = \X_e(m,n;r,t)$, where $m,n \geq 4$ are even integers and $r\in \ZZ_n^*$, 
$t \in \ZZ_n$ are such that (\ref{eq:Xparamcond}) holds. Let $G$ be
a subgroup of automorphisms of $X$ such that $H \leq G \leq \Aut X$. Let $\C$ be a union of $G$-orbits
of cycles of $X$ and let $P$ and $Q$ be any two paths of $X$ permutable by some element of $G$. Then
$v(\C,P) = v(\C,Q)$.
\end{proposition}

\begin{proof}
Let $Bip$ be the bipartite graph having as one bipartition set the set of cycles $\C$ and the
other bipartition set the $G$-orbit $\P$ of $P$, such that a cycle $C \in \C$ is adjacent to
a path $R \in \P$ if $C$ contains $R$ as a subgraph. Since $\C$ is a union of $G$-orbits of cycles,
there is a natural action of $G$ on $Bip$. Clearly for any $R \in \P$ the frequency $v(\C,R)$ is precisely
the valency of $R$ in $Bip$. Thus since $\P$ is a $G$-orbit containing $P$ and $Q$, the result follows.
\end{proof}\bigskip

Let $\C$ be a union of $H$-orbits of cycles of $X$. Let $x$ be any of the symbols $\{a^+, a^-, g, z\}$.
By Proposition~\ref{pro:equalfreq} we can now define $v_x(\C)$ to be the frequency $v(\C,P)$ where
$P$ is any $2$-path of type $x$. Note that by Proposition~\ref{pro:anchors} we have $v_{a^+}(\C) = v_{a^-}(\C)$,
so we set $v_a(\C) = v_{a^+}(\C) = v_{a^-}(\C)$. In accordance
with the above notation we define the frequencies $v_{a}(d)$, $v_g(d)$ and $v_z(d)$ to be the
respective $d$-frequencies of anchors, glides and zigzags, respectively. The next lemma gives an
easy way of calculating these frequencies and we will be using it throughout the rest of the paper
without special reference to it.

\begin{lemma}
\label{le:calcfreq}
Let $X = \X_e(m,n;r,t)$, where $m,n \geq 4$ are even integers and $r\in \ZZ_n^*$, 
$t \in \ZZ_n$ are such that (\ref{eq:Xparamcond}) holds. Let $\C$ be an
$H$-orbit of cycles of $X$ and let $x$ be any symbol from the set $\{a,g,z\}$. Let $C \in \C$ and suppose
$C$ contains $\varepsilon_{x,C}$ $2$-paths of type $x$. Then
$$ v_x(\C) = \frac{|\C|\cdot \varepsilon_{x,C}}{2mn}. $$
\end{lemma}

\begin{proof}
Let $P$ be any $2$-path of type $x$ and let $\P$ be its $H$-orbit. We let $Bip$ be the bipartite
graph having as the two bipartition sets $\C$ and $\P$ such that $C' \in \C$ is adjacent to $P'\in \P$
whenever $C'$ contains $P'$ as a subgraph. Since $H$ acts transitively on each of the two bipartition sets, 
using a simple
counting argument and Proposition~\ref{pro:equalfreq}, we can 
show that $|\C|\cdot \varepsilon_{x,C} = |\P|\cdot v_x(\C)$ if $x \in \{g,z\}$. 
By $(\ref{eq:orblen})$, the result follows.
Suppose now that $x = a$. Then $C$ contains $\frac{\varepsilon_{a,C}}{2}$ positive
and $\frac{\varepsilon_{a,C}}{2}$ negative anchors by Proposition~\ref{pro:anchors}. 
By the above argument we have
$|\C|\cdot \frac{\varepsilon_{a,C}}{2} = |Anc^+ X|\cdot v_{a^+}(\C) = |Anc^+ X|\cdot v_{a}(\C)$ which,
in view of $(\ref{eq:orblen})$, completes the proof.
\end{proof}\bigskip

We end this section with two propositions 
that will be used in our investigation of the arc-transitivity and half-arc-transitivity of $X$ later on.

\begin{proposition}
\label{pro:halftrcond}
Let $X = \X_e(m,n;r,t)$, where $m,n \geq 4$ are even integers and $r\in \ZZ_n^*$, 
$t \in \ZZ_n$ are such that  (\ref{eq:Xparamcond}) holds, and furthermore  $r^2 \neq \pm 1$. 
Let $\C$ be an $\Aut X$-orbit of $d$-cycles of $X$ for some $d \geq 4$. 
If $v_a(\C) \neq v_g(\C)$ and $v_a(\C) \neq v_z(\C)$ then $X$ is half-arc-transitive. 
In particular, if $v_a(d) \neq v_g(d)$ and $v_a(d) \neq v_z(d)$
for some $d \geq 4$, then $X$ is half-arc-transitive.
\end{proposition}

\begin{proof}
Suppose on the contrary that $X$ is arc-transitive. By Proposition~\ref{pro:arc-tr} some anchor of
$X$ is mapped by some automorphism of $X$ to a glide or to a zigzag. But then Proposition~\ref{pro:equalfreq}
implies that either $v_a(\C) = v_g(\C)$ or $v_a(\C) = v_z(\C)$, a contradiction.
\end{proof}

\begin{proposition}
\label{pro:theisomorphisms}
Let $X = \X_e(m,n;r,t)$, where $m,n \geq 4$ are even integers and $r\in \ZZ_n^*$, 
$t \in \ZZ_n$ are such that (\ref{eq:Xparamcond}) holds.
Let $Y = \X_e(m,n;-r,t+r+r^3+\cdots + r^{m-1})$ and let $Z = \X_e(m,n;r^{-1},t)$.
Then $X \iso Y \iso Z$.
Moreover, there exist isomorphisms of graphs
$\varphi : X \to Y$ and $\psi : X \to Z$, such that $Anc X\varphi = Anc Y$, 
$Gli X\varphi = Zig Y$ and $Zig X\varphi = Gli Y$, and 
$Anc X \psi = Anc Z$, $Gli X\psi = Gli Z$ and $Zig X \psi = Zig Z$.
\end{proposition}

\begin{proof}
Let us first show that $X \iso Z$. Denote the
vertex set of $Z$ by $\{v_i^j\ |\ i \in \ZZ_m,\ j \in \ZZ_n\}$ with edges as usual.
Let $\psi : X \to Z$ be the mapping defined by the rule
$$ u_i^j\psi = \left\{\begin{array}{lll}
		v_0^{-rj}	& ;	& i = 0,\ j \in \ZZ_n \\
		v_{m-i}^{-rj-t}	& ; 	& i \in \ZZ_m \setminus \{0\},\ j \in \ZZ_n .
		\end{array} \right.$$
Clearly $\psi$ is bijective. Since $r^m = 1$, we have $r^{-1} = r^{m-1}$, and so
$r^{m-i} = (r^{-1})^{i}$. Using $rt = t$ from (\ref{eq:Xparamcond}), it is 
straightforward to check that $\psi$ is a graph homomorphism. It is also clear that
$Anc X \psi = Anc Z$, $Gli X\psi = Gli Z$ and $Zig X \psi = Zig Z$.

Let us now construct an isomorphism $\varphi : X \to Y$, 
so that the required conditions are met. Denote the vertex set of $Y$ by 
$\{v_i^j\ |\ i \in \ZZ_m,\ j \in \ZZ_n\}$ with edges as usual. Let $\varphi : X \to Y$
be the mapping defined by the rule
$$ u_i^j\varphi = \left\{\begin{array}{lll}
		v_i^j				& ;	& i \in \{0,1\},\ j \in \ZZ_n\\
		v_i^{j-r-r^3-\cdots - r^{i-1}}	& ; 	& i \in \ZZ_m \setminus \{0\},\ i \equiv 0 \pmod 2,\ j \in \ZZ_n\\
		v_i^{j-r-r^3-\cdots - r^{i-2}}	& ; 	& i \in \ZZ_m \setminus \{1\},\ i \equiv 1 \pmod 2,\ j \in \ZZ_n .
		\end{array}\right. $$
Clearly $\varphi$ is a bijection. We now show that $\varphi$ is a graph homomorphism.
For instance, every edge connecting a vertex of $X_1$ to a vertex of $X_2$ (recall that $X_i = \{u_i^j\ |\ j \in \ZZ_n\}$) 
is of the form $u_1^ju_2^j$ or $u_1^ju_2^{j+r}$. Since $u_1^j\varphi = v_1^j$, 
$u_2^j\varphi = v_2^{j-r}$ and $u_2^{j+r}\varphi = v_2^{j}$, the two edges in question are mapped to edges of $Y$.
Similarly an edge connecting a vertex of $X_{m-1}$ to a vertex of $X_0$ is of the form
$u_{m-1}^ju_0^{j+t}$ or $u_{m-1}^ju_0^{j+r^{m-1}+t}$. Since $m-1 \equiv 1 \pmod 2$, 
we have $u_{m-1}^j\varphi = v_{m-1}^{j-r-r^3-\cdots - r^{m-3}}$ and
$u_0^{j+t} = v_0^{j+t}$, $u_0^{j+r^{m-1}+t}\varphi = v_0^{j+r^{m-1} + t}$.
But since 
\begin{equation}\label{equ:conditionD}
j-r-r^3-\cdots - r^{m-3} + (-r)^{m-1} + t + r + r^3 + \cdots + r^{m-1} = j + t
\end{equation}
we have  
$j-r-r^3-\cdots - r^{m-3} + t + r + r^3 + \cdots + r^{m-1} = j + r^{m-1} + t$, and so the two edges in question
are mapped to edges of $Y$.
We now only need to see that the edges connecting the vertices of $X_i$ to those of $X_{i+1}$, 
where $2 \leq i \leq m-2$, are mapped to edges of $Y$. We leave this
to the reader. 
Let $Y_i = \{v_i^j\ |\ j \in \ZZ_n\}$. Clearly $\varphi$ maps each 
$X_i$, $i \in \ZZ_m$, to $Y_i$, and so $\varphi$ maps $Anc X$ to $Anc Y$.
To show that $\varphi$ maps $Gli X$ to $Zig Y$ it thus suffices to show
that for every vertex $u$ of $X$,  $\varphi$ maps a glide with $u$ as its internal vertex to a zigzag of $Y$.
Let $j \in \ZZ_n$ be arbitrary. Then the glide $u_0^ju_1^ju_2^j$ is mapped to the zigzag $v_0^jv_1^jv_2^{j-r}$, 
the glide $u_1^ju_2^ju_3^j$ is mapped to the zigzag $v_1^jv_2^{j-r}v_3^{j-r}$, and so on.
Finally, the glide $u_{m-2}^ju_{m-1}^ju_0^{j+t}$ is mapped to the $2$-path 
$v_{m-2}^{j-r-r^3-\cdots -r^{m-3}}v_{m-1}^{j-r-r^3-\cdots -r^{m-3}}v_0^{j+t}$ of $Y$ which in view of
(\ref{equ:conditionD})
is a zigzag . Similarly the glide $u_{m-1}^ju_0^{j+t}u_1^{j+t}$ is mapped to a zigzag 
$v_{m-1}^{j-r-r^3-\cdots -r^{m-3}}v_0^{j+t}v_1^{j+t}$ of $Y$. 
Thus $\varphi$ maps  $Gli X$ to  $Zig Y$, and hence it also maps $Zig X$ to  $Gli Y$.
\end{proof}


\section{The $8$-cycles}
\label{sec:8cycles}

\indent
Throughout this section we let $m,n,r,t,X,H,\rho,\sigma$ and $\tau$ be as in Section~\ref{sec:notation}.
Since our goal is to determine whether or not $X$ is half-arc-transitive, we can (in view of 
Proposition~\ref{pro:arc-tr}) assume that $r^2 \neq \pm 1$. Since $n$ is even, 
$r \in \ZZ_n^*$ implies that $r$ is odd. Suppose now that $2r = \pm 2$. Then
$2(r \mp 1) = 0$ implies $r \mp 1 = \frac{n}{2}$ (or else $r^2 = 1$), which must be even, and so $n \equiv 0 \pmod 4$.
Thus $r^2 = ((r \mp 1) \pm 1)^2 = (\frac{n}{2})^2 \pm n + 1 = 1$, which was assumed not to be true. Therefore
\begin{equation}
\label{eq:2r}
	2r \neq \pm 2.
\end{equation}
Since $r \in \ZZ_n^*$ and $r^2 \neq \pm 1$ we can also assume
\begin{equation}
\label{eq:n}
	n \geq 14.
\end{equation}
Throughout the rest of the paper we will constantly be relying on these two observations.

We now determine the set of all possible $8$-cycles of $X$. Recall that 
$X_i = \{u_i^j\ |\ j \in \ZZ_n\}$, where $i \in \ZZ_m$, are the $H$-attachment sets of $X$
as well as the orbits of $\rho$. For brevity reasons we will simply call them orbits in the discussion bellow.
We distinguish different cases depending on
the number of orbits the $8$-cycle in question lies on and determine the
possible traces such $8$-cycles can have. 

\begin{itemize}
\item No $8$-cycle of $X$ lies on two consecutive
orbits, for then $n = 4$, which contradicts (\ref{eq:n}). 

\item The trace of an $8$-cycle on three consecutive
orbits could either be $anananan$, $a^5nan$ or $a^3na^3n$. 
It is easy to see, however, that an $8$-cycle of trace $a^3na^3n$
would exist if and only if either $2r+2 = 0$ or $2r-2 = 0$, which cannot occur by (\ref{eq:2r}). 
Therefore, the $8$-cycles on three consecutive orbits are either of trace $anananan$,
or of trace $a^5nan$.

\item Let us consider the possible traces of $8$-cycles lying on four orbits. 
Suppose first that $m > 4$. It is easy to see that in this case the trace of an $8$-cycle 
is either $a^2nan^2an$ or $a^3n^2an^2$. Suppose now that $m = 4$. Then a careful examination
of possible cases shows, that in addition to the above two traces an $8$-cycle can also have 
trace $n^8$, $anan^5$, $a^4n^4$, $a^2na^2n^3$ or $a^2n^2a^2n^2$.

\item Clearly the trace of an $8$-cycle on five orbits is $an^3an^3$.

\item An $8$-cycle lying on six consecutive orbits can exist only when $m = 6$ in which case its trace is $a^2n^6$.

\item Finally, an $8$-cycle on $8$ orbits can only be a coiled one and thus of trace $n^8$.
Clearly $m = 8$ in this case.
\end{itemize}

We now determine all possible $8$-cycles of traces $anananan$, $a^5nan$, $a^3n^2an^2$, $a^2nan^2an$
and $an^3an^3$. As the other $8$-cycles can only exist in special cases when $m = 4, 6$ or $8$, we
will deal with them in subsequent sections, where these special cases of $m$ are considered. 
For each of the possible traces a careful examination of the
possible cycles needs to be undertaken. To indicate how this is done, we consider the trace $anananan$
and leave the rest to the reader. The possible $H$-orbits of $8$-cycles are collected in Table~\ref{tab:8cycles}
bellow. Each row corresponds to one $H$-orbit of $8$-cycles. We include the following information: the trace 
and the code of the $8$-cycles of the $H$-orbit in question, 
its representative, the necessary and sufficient condition under 
which the $H$-orbit exists and the length of the $H$-orbit. The reader should note that in case
$m = 4$ the orbit $X_4$ is in fact equal to $X_0$, so 
some of the proposed $H$-orbits of $8$-cycles of trace $an^3an^3$ might not exist
even if the condition stated in the table holds. Namely, the 
$8$ vertices of the given representative might not be distinct 
(this can only occur when $m = 4$ and the girth of $X$ is $4$). 
The additional condition for the existence
of these $H$-orbits in case $m = 4$ is thus that the vertex of form $u_4^j$, which in this case is
actually $u_0^{j+t}$,
is different from $u_0^0$. It turns out that this never causes problems in our investigation though, 
so we do not state this conditions in the table to improve readability.

Let us now investigate the possible $H$-orbits of $8$-cycles of trace $anananan$.
Consider an $8$-cycle $C$ containing the negative anchor $u_1^1u_0^0u_1^0$. 
If both $u_1^0$ and $u_1^1$ correspond to glides on $C$, then $C$ contains the path
$u_1^{1+r}u_2^{1+r}u_1^1u_0^0u_1^0u_2^0u_1^{-r}$, so since $C$ is an $8$-cycle either $1+r+1 = -r$,
which contradicts (\ref{eq:2r}), or $-r+1 = 1+r$, and so $2r = 0$, which is also impossible.
Suppose now that the vertices $u_1^0$ and $u_1^1$ correspond to zigzags of $C$. Then $C$ contains the path
$u_1^{1-r}u_2^1u_1^1u_0^0u_1^0u_2^ru_1^r$, and so either $1-r+1 = r$ or $r+1 = 1-r$, which are both impossible.
Suppose finally that the vertices $u_1^0$ and $u_1^1$ correspond to one glide and one zigzag of $C$. 
With no loss of generality (since $\tau \in H$) we can assume that $u_1^0$ corresponds to a glide. 
Thus $C$ contains the path
$u_1^{1-r}u_2^1u_1^1u_0^0u_1^0u_2^0u_1^{-r}$, and so the only possibility for $C$
to be a cycle is that the remaining vertex is $u_0^{-r}$. Note that this $8$-cycle always exists in $X$. 
Moreover, the code of $C$ is $agazagaz$ and
the $H$-orbit of $C$ has length $mn$, since the only automorphism of $H$ that
fixes $C$ setwise is $\tau\rho^{-r}$. 
We call the $8$-cycles of this $H$-orbit the {\em generic} $8$-cycles of $X$.
They correspond to row 1 of Table~\ref{tab:8cycles}. The investigation of possible $8$-cycles
of other traces is done in a similar manner and is left to the reader.

\begin{table}[!tb]
\begin{footnotesize}
\begin{center}
\begin{tabular}{@{}c|c|l|c|c|c@{}}
	Row	& Trace		& A representative 	& Code	& Condition	& Orbit length \\
	\hline  & & & & & \\
   $1$ & $anananan$ & $u_0^0u_1^0u_2^0u_1^{-r}u_0^{-r}u_1^{1-r}u_2^{1}u_1^1$ & $agazagaz$ & none & $mn$ \\ 
   & & & & & \\
   $2$ & $a^5nan$   & $u_0^0u_1^1u_0^1u_1^2u_0^2u_1^3u_2^3u_1^{3-r}$	& $a^5zaz$ & $3-r=0$ & $mn$ \\
   $3$ & $a^5nan$   & $u_0^0u_1^1u_0^1u_1^2u_0^2u_1^3u_2^{3+r}u_1^{3+r}$	& $a^5gag$ & $3+r=0$ & $mn$ \\
   $4$ & $a^5nan$   & $u_0^0u_1^0u_2^ru_1^ru_2^{2r}u_1^{2r}u_2^{3r}u_1^{3r}$ & $a^5zaz$ & $1-3r = 0$ & $mn$\\
   $5$ & $a^5nan$   & $u_0^0u_1^1u_2^{1+r}u_1^{1+r}u_2^{1+2r}u_1^{1+2r}u_2^{1+3r}u_1^{1+3r}$ & $a^5gag$ & $1+3r = 0$ & $mn$\\    
   & & & & & \\
   $6$ & $a^2nan^2an$ & $u_0^0u_1^0u_2^{0}u_1^{-r}u_2^{-r}u_3^{-r}u_2^{-r-r^2}u_1^{-2r-r^2}$ & $a^2gag^2ag$ & $1+2r+r^2= 0$ & $2mn$\\
   $7$ & $a^2nan^2an$ & $u_0^0u_1^0u_2^{0}u_1^{-r}u_2^{-r}u_3^{-r+r^2}u_2^{-r+r^2}u_1^{-2r+r^2}$ & $a^2zazgag$ & $1+2r-r^2= 0$ & $2mn$\\
   $8$ & $a^2nan^2an$ & $u_0^0u_1^0u_2^{r}u_1^{r}u_2^{2r}u_3^{2r}u_2^{2r-r^2}u_1^{2r-r^2}$ & $a^2zaz^2az$ & $1-2r+r^2= 0$ & $2mn$\\
   $9$ & $a^2nan^2an$ & $u_0^0u_1^0u_2^{r}u_1^{r}u_2^{2r}u_3^{2r+r^2}u_2^{2r+r^2}u_1^{2r+r^2}$ & $a^2gagzaz$ & $1-2r-r^2= 0$ & $2mn$\\
   & & & & & \\
   $10$& $a^3n^2an^2$ & $u_0^0u_1^1u_0^1u_1^2u_2^2u_3^2u_2^{2-r^2}u_1^{2-r-r^2}$ & $a^3zgagz$ & $2-r-r^2=0$ & $mn$\\
   $11$& $a^3n^2an^2$ & $u_0^0u_1^1u_0^1u_1^2u_2^2u_3^{2+r^2}u_2^{2+r^2}u_1^{2-r+r^2}$ & $a^3z^2az^2$ & $2-r+r^2=0$ & $mn$\\
   $12$& $a^3n^2an^2$ & $u_0^0u_1^1u_0^1u_1^2u_2^{2+r}u_3^{2+r}u_2^{2+r-r^2}u_1^{2+r-r^2}$ & $a^3gzazg$ & $2+r-r^2=0$ & $mn$\\
   $13$& $a^3n^2an^2$ & $u_0^0u_1^1u_0^1u_1^2u_2^{2+r}u_3^{2+r+r^2}u_2^{2+r+r^2}u_1^{2+r+r^2}$ & $a^3g^2ag^2$ & $2+r+r^2=0$ & $mn$\\
   & & & & & \\
   $14$& $a^3n^2an^2$ & $u_0^0u_1^0u_2^0u_3^0u_2^{-r^2}u_3^{-r^2}u_2^{-2r^2}u_1^{-r-2r^2}$ & $a^3g^2ag^2$ & $1+r+2r^2=0$ & $mn$\\
   $15$& $a^3n^2an^2$ & $u_0^0u_1^1u_2^1u_3^1u_2^{1-r^2}u_3^{1-r^2}u_2^{1-2r^2}u_1^{1-r-2r^2}$ & $a^3gzazg$ & $1-r-2r^2=0$ & $mn$\\
   $16$& $a^3n^2an^2$ & $u_0^0u_1^0u_2^ru_3^ru_2^{r-r^2}u_3^{r-r^2}u_2^{r-2r^2}u_1^{r-2r^2}$ & $a^3z^2az^2$ & $1-r+2r^2=0$ & $mn$\\
   $17$& $a^3n^2an^2$ & $u_0^0u_1^1u_2^{1+r}u_3^{1+r}u_2^{1+r-r^2}u_3^{1+r-r^2}u_2^{1+r-2r^2}u_1^{1+r-2r^2}$ & $a^3zgagz$ & $1+r-2r^2=0$ & $mn$\\
   & & & & & \\
   $18$& $an^3an^3$ & $u_0^0u_1^0u_2^0u_3^0u_4^{r^3}u_3^{r^3}u_2^{r^3}u_1^{r^3}$ & $ag^2zag^2z$ & $1-r^3=0$ & $2mn$\\
   $19$& $an^3an^3$ & $u_0^0u_1^0u_2^0u_3^{r^2}u_4^{r^2+r^3}u_3^{r^2+r^3}u_2^{r^3}u_1^{r^3}$ & $agzgaz^3$ & $1-r^3=0$ & $2mn$\\
   & & & & & \\
   $20$& $an^3an^3$ & $u_0^0u_1^0u_2^0u_3^0u_4^{0}u_3^{-r^3}u_2^{-r^3}u_1^{-r^3}$ & $ag^3azgz$ & $1+r^3=0$ & $2mn$\\
   $21$& $an^3an^3$ & $u_0^0u_1^0u_2^0u_3^{r^2}u_4^{r^2}u_3^{r^2-r^3}u_2^{-r^3}u_1^{-r^3}$ & $agz^2agz^2$ & $1+r^3=0$ & $2mn$\\
   & & & & & \\
   $22$& $an^3an^3$ & $u_0^0u_1^0u_2^0u_3^0u_4^{0}u_3^{-r^3}u_2^{-r^2-r^3}u_1^{-r-r^2-r^3}$ & $ag^3ag^3$ & $1+r+r^2+r^3=0$ & $mn$\\
   $23$& $an^3an^3$ & $u_0^0u_1^0u_2^0u_3^0u_4^{r^3}u_3^{r^3}u_2^{-r^2+r^3}u_1^{-r-r^2+r^3}$ & $ag^2zazg^2$ & $1+r+r^2-r^3=0$ & $mn$\\
   $24$& $an^3an^3$ & $u_0^0u_1^0u_2^0u_3^{r^2}u_4^{r^2}u_3^{r^2-r^3}u_2^{r^2-r^3}u_1^{-r+r^2-r^3}$ & $agz^2az^2g$ & $1+r-r^2+r^3=0$ & $mn$\\
   $25$& $an^3an^3$ & $u_0^0u_1^0u_2^0u_3^{r^2}u_4^{r^2+r^3}u_3^{r^2+r^3}u_2^{r^2+r^3}u_1^{-r+r^2+r^3}$ & $agzgagzg$ & $1+r-r^2-r^3=0$ & $mn$\\  
   $26$& $an^3an^3$ & $u_0^0u_1^0u_2^ru_3^ru_4^{r}u_3^{r-r^3}u_2^{r-r^2-r^3}u_1^{r-r^2-r^3}$ & $az^2gagz^2$ & $1-r+r^2+r^3=0$ & $mn$\\
   $27$& $an^3an^3$ & $u_0^0u_1^0u_2^ru_3^ru_4^{r+r^3}u_3^{r+r^3}u_2^{r-r^2+r^3}u_1^{r-r^2+r^3}$ & $az^3az^3$ & $1-r+r^2-r^3=0$ & $mn$\\
   $28$& $an^3an^3$ & $u_0^0u_1^0u_2^ru_3^{r+r^2}u_4^{r+r^2}u_3^{r+r^2-r^3}u_2^{r+r^2-r^3}u_1^{r+r^2-r^3}$ & $azgzazgz$ & $1-r-r^2+r^3=0$ & $mn$\\
   $29$& $an^3an^3$ & $u_0^0u_1^0u_2^ru_3^{r+r^2}u_4^{r+r^2+r^3}u_3^{r+r^2+r^3}u_2^{r+r^2+r^3}u_1^{r+r^2+r^3}$ & $azg^2ag^2z$ & $1-r-r^2-r^3=0$ & $mn$\\
\end{tabular}
\caption{$H$-orbits of $8$-cycles of certain traces.}
\label{tab:8cycles}
\end{center}
\end{footnotesize}
\end{table}	

The next lemmas and a corollary are the first steps towards determining the actual 
$H$-orbits of $8$-cycles of $X$.

\begin{lemma}
\label{le:possHorbits}
Let $X = \X_e(m,n;r,t)$, where $m,n \geq 4$ are even integers and $r\in \ZZ_n^*$, 
$t \in \ZZ_n$ are such that (\ref{eq:Xparamcond}) holds, and furthermore $r^2 \neq \pm 1$. 
Then $X$ can have at most
\begin{itemize}
	\item[(i)] one $H$-orbit of $8$-cycles of trace $a^5nan$;
	\item[(ii)] one $H$-orbit of $8$-cycles of trace $a^2nan^2an$;
	\item[(iii)] one $H$-orbit of $8$-cycles of trace $a^3n^2an^2$ from rows $10$ to $13$ of 
	Table~\ref{tab:8cycles};
	\item[(iv)] one $H$-orbit of $8$-cycles of trace $a^3n^2an^2$ from rows $14$ to $17$ of 
	Table~\ref{tab:8cycles};
	\item[(v)] one $H$-orbit of $8$-cycles of trace $an^3an^3$ from rows $22$ to $29$ of 
	Table~\ref{tab:8cycles} unless it has two such orbits in which case either 
	$r^2 = \frac{n}{2} - 1$ and the two $H$-orbits correspond to rows $22$ and $27$, or
	$r^2 = \frac{n}{2} + 1$ and the two $H$-orbits correspond to rows $25$ and $28$;
	\item[(vi)] three $H$-orbits of $8$-cycles of trace $an^3an^3$ and if it does have three
	such $H$-orbits, then $r^3 = \pm 1$.
\end{itemize}
\end{lemma}

\begin{proof}
Using (\ref{eq:n}) and the fact that $r \in \ZZ_n^*$, the claims (i) through (iv) are easily verified. Suppose now
that the $H$-orbit of row $22$ exists. Using (\ref{eq:2r}) it can be seen that no additional condition of rows
$22$-$29$, except for that of row $27$, can hold.
If the condition of row $27$ does hold, then $2(1+r^2)=0$ and since 
$r^2 \neq -1$, we have $1+r^2 = \frac{n}{2}$, as claimed. Furthermore, if the $H$-orbit of
row $23$ exists, then the only other $H$-orbit of rows $22$-$29$ that could exist is the one of row $26$. But if
that was true, then $2(1+r^2)=0$ and $2(r-r^3)=0$, and so $2(1-r^2)=0$, which implies $4=0$, contradicting $(\ref{eq:n})$. 
A similar argument shows that the $H$-orbit of row $24$ cannot exist simultaneously with 
any other $H$-orbit of rows $22$-$29$.
Finally, if the $H$-orbit of row $25$ exists, then the only other $H$-orbit of rows $22$-$29$ that
can exist is the one of row $28$. In this case $2(1-r^2)=0$ and thus $r^2 \neq 1$ implies that $1-r^2 = \frac{n}{2}$,
proving (v). As for (vi), it is clear that we cannot have $r^3 = 1$ and $r^3 = -1$. By (v) we can have
two $H$-orbits corresponding to rows $22$-$29$ only if $r^2 = \frac{n}{2} \pm 1$. But then
$r^3 = r(\frac{n}{2} \pm 1) = \frac{n}{2} \pm r \neq \pm 1$, or else $2r = \pm 2$, which contradicts (\ref{eq:2r}).
\end{proof}\bigskip

Using Lemma~\ref{le:calcfreq} the next corollary is now straightforward. 

\begin{corollary}
\label{cor:summs}
Let $X = \X_e(m,n;r,t)$, where $m,n \geq 4$ are even integers and $r\in \ZZ_n^*$, 
$t \in \ZZ_n$ are such that  (\ref{eq:Xparamcond}) holds, and furthermore  $r^2 \neq \pm 1$. Let
$\C_1$ be a union of $H$-orbits of $8$-cycles of $X$ whose traces are either $a^2nan^2an$ or $a^3n^2an^2$ and let
$\C_2$ be a union of $H$-orbits of $8$-cycles of $X$ of trace $an^3an^3$. Then $v_a(\C_1) = v_g(\C_1) + v_z(\C_1)$
is an even number and $v_a(\C_2) = \frac{1}{3}(v_g(\C_2) + v_z(\C_2)) \leq 5$.
\end{corollary}

\begin{lemma}
\label{le:C1orb}
Let $X = \X_e(m,n;r,t)$, where $m,n \geq 4$ are even integers and $r\in \ZZ_n^*$, 
$t \in \ZZ_n$ are such that  (\ref{eq:Xparamcond}) holds, and furthermore  $r^2 \neq \pm 1$.
If there exist at least two $H$-orbits of $8$-cycles of $X$ corresponding to rows $6$-$17$
of Table~\ref{tab:8cycles}, then there exists an $8$-cycle of trace $a^5nan$ in $X$ or up to
isomorphisms given by Proposition~\ref{pro:theisomorphisms} one of the following holds:
\begin{itemize}
	\item[(i)] $r^3 = \pm 1$;
	\item[(ii)] $r = 5$ and $n = 22$, and so $r$ is of order $5$ in $\ZZ_n^*$;
	\item[(iii)] $r = 5$ and $n = 28$, and so $r$ is of order $6$ in $\ZZ_n^*$.
\end{itemize}
\end{lemma}

\begin{proof}
In view of the isomorphism $X \iso  \X_e(m,n;r^{-1},t)$, given by Proposition~\ref{pro:theisomorphisms},
Lemma~\ref{le:possHorbits} implies, that we can assume that an $H$-orbit corresponding to one of the rows
$10$-$13$ of Table~\ref{tab:8cycles} exists. 

Suppose first that some $H$-orbit of $8$-cycles of trace $a^2nan^2an$ also exists.
In view of the isomorphism $X \iso \X_e(m,n;-r,t+r+r^3+\cdots + r^{m-1})$, we can assume that
an $H$-orbit corresponding to row $6$ or row $7$ exists. Suppose first it is the
one of row $6$. Then the condition of row $10$ implies $3+r=0$, that of row $11$ implies $1-3r=0$, so
an $8$-cycle of trace $a^5nan$ exists in these two cases. The condition of row 
$12$ implies $3+3r=0$, and since also $(1+r)^2=0$, we
get $r^3 = (r+1-1)^3 = -1$. Finally, the condition of row $13$ implies $1-r=0$, 
which is impossible. Suppose now that the $H$-orbit of row $7$ exists. Then the conditions
of rows $10$ and $11$ give $1-3r=0$ and $3+r=0$, respectively, whereas the one of row~$12$ gives $1-r=0$,
which contradicts $r^2 \neq 1$. 
The condition of row~$13$, however, implies $3+3r=0$, and so $0 = 3(2+r+r^2) = 6 + r(3+3r) = 6$,
which contradicts $(\ref{eq:n})$.

Suppose now that an $H$-orbit corresponding to rows $14$-$17$ of Table~\ref{tab:8cycles} exists (recall
that we are already assuming that an $H$-orbit corresponding to one of the rows $10$-$13$ exists).
In view of Proposition~\ref{pro:theisomorphisms} we can assume that either
the $H$-orbit corresponding to row $10$ or the one corresponding to row $11$ exists.
Suppose first it is the one of row $10$. Then the condition of row $14$ implies $r=5$, and so 
$28 \equiv 0 \pmod n$. By $(\ref{eq:n})$, $n$ is either $14$ or $28$. If $n = 14$ then $1-3r=0$ and
if $n = 28$ the order of $r$ in $\ZZ_n^*$ is $6$. The condition of
row $15$ gives $3-r=0$. The one of row $16$ gives $5-3r=0$. Thus $0 = 6-3r-3r^2 = -4-2r$, which implies
$r=9$, and so $n=22$. Note that $r^{-1} = 5$ in this case. The condition of row~$17$,
on the other hand, gives $3-3r=0$ and $1-2r+r^2=0$. Thus $r^3 = (r-1+1)^3 = 1$ in this case.
Suppose finally that the $H$-orbit corresponding to row~$11$ of Table~\ref{tab:8cycles} exists.
Then the condition of row~$14$ gives $3-3r=0$, and so $0=3+3r+6r^2=12$, contradicting $(\ref{eq:n})$.
Furthermore, the condition of row~$15$ implies $5-3r=0$, and so $0 = 6-3r+3r^2 = 6+2r$. Thus $r = 11$,
forcing $28 \equiv 0 \pmod n$. Therefore, 
$n$ is either $14$ or $28$. In the first case $3+r=0$ and in the second case
$-r^{-1} = 5$. 
Finally, the condition of row~$16$ implies $3-r=0$ and that of row~$17$ implies $r=5$, forcing $n = 22$.
\end{proof}\bigskip

The following lemma, whose proof depends heavily on 
\cite[Lemma~2.1.]{DM05}, will play an important role in the investigation
of half-arc-transitivity of the graphs $\X_e(m,n;r,t)$ in the subsequent sections.

\begin{lemma}
\label{le:quartic}
Let $X = \X_e(m,n;r,t)$, where $m,n \geq 4$ are even integers and $r\in \ZZ_n^*$, 
$t \in \ZZ_n$ are such that  (\ref{eq:Xparamcond}) holds, and furthermore  $r^2 \neq \pm 1$. Suppose 
$X$ is arc-transitive. Then the $\Aut X$-orbit of the generic
$8$-cycles contains $8$-cycles with two consecutive anchors.
\end{lemma}

\begin{proof}
Let $D(X)$ be one of the two oriented graphs corresponding to
the half-arc-transitive action of $H$ on $X$.
Proposition~\ref{pro:arc-tr} implies that the automorphism group of $X$ does not fix the set of
anchors of $D(X)$. 
Thus there exists an automorphism $\varphi$ of $X$, mapping some anchor either to a glide or to a zigzag.
Therefore $\varphi$ neither preserves nor reverses the orientation of every edge of $D(X)$. 
For the purposes of this proof only we let the refinement of a trace of a path $W$ be the sequence obtained
from its trace by replacing each $a$ by $a^+$ or $a^-$, depending on 
what type of anchor the symbol $a$ corresponds to. The possible
refinements of traces of $3$-paths are thus $nn$, $a^+a^-$, $na^+$ and $na^-$.
By \cite[Lemma~2.1.]{DM05}, it is either true 
that every automorphism of $X$ preserves the orientation
of every edge of $D(X)$ or reverses the orientation of every edge of $D(X)$,
or that for any two of the four possible refinements of traces $t_1$ and $t_2$ of $3$-paths of $X$ there exist 
$3$-paths $P_1$ and $P_2$ with refinements of traces $t_1$ and $t_2$, respectively, 
and an automorphism of $X$ mapping $P_1$ to $P_2$. 
Therefore, since the former is not true in our case, there exist $3$-paths $P_1$
and $P_2$, such that the refinement of the trace of $P_1$ is $na^+$ and 
the refinement of the trace of $P_2$ is $a^+a^-$, and there exists an automorphism 
$\varphi$ of $X$ mapping $P_1$ to $P_2$. But there are only two $H$-orbits of $3$-paths having
refinements of traces $na^+$, one whose $3$-paths have refinements of codes $ga^+$ and the 
other whose $3$-paths have refinements of codes $za^+$. 
Since a $3$-path of each of these two refinements of codes lies on some generic $8$-cycle, 
the $3$-path $P_1$ lies on some generic $8$-cycle, and the result follows.
\end{proof}


\section{The general case}

\indent
Throughout this section we let $m,n,r,t,X,H,\rho,\sigma$ and $\tau$ be as in Section~\ref{sec:notation}.
Recall that this implies that the conditions (\ref{eq:Xparamcond}) are satisfied. In view of Proposition~\ref{pro:r^2}
we can assume $r^2 \neq \pm 1$.
As already mentioned in the previous section an $8$-cycle of trace 
$n^8$, $a^2n^6$, $anan^5$, $a^4n^4$, $a^2na^2n^3$ or $a^2n^2a^2n^2$ can only exist if $m$ is one of $4,6$ or $8$.
It thus seems only natural to first consider graphs that do not have $8$-cycles of the above traces.
It is the aim of this section to show (see Lemma~\ref{le:general} bellow) that the graph $X$
is half-arc-transitive in this case.  

\begin{lemma}
\label{le:genequal}
Let $X = \X_e(m,n;r,t)$, where $m,n \geq 4$ are even integers and $r\in \ZZ_n^*$, 
$t \in \ZZ_n$ are such that  (\ref{eq:Xparamcond}) holds, and furthermore  $r^2 \neq \pm 1$.
If $X$ contains no $8$-cycle of trace $n^8$, $a^2n^6$, $anan^5$, $a^4n^4$, $a^2na^2n^3$ or $a^2n^2a^2n^2$
and $v_a(8) = v_g(8) = v_z(8)$, then $X$ is half-arc-transitive.
\end{lemma}

\begin{proof}
Note that the assumptions imply that all the possible $H$-orbits of $8$-cycles of $X$ are listed in 
Table~\ref{tab:8cycles}.
Let $\C$ denote the set of all $8$-cycles of $X$, 
let $\C_1$ denote the set of $8$-cycles of $X$ of traces $a^2nan^2an$ and $a^3n^2an^2$ and let $\C_2$
denote the set of $8$-cycles of $X$ of trace $an^3an^3$. We claim that no $8$-cycle of trace $a^5nan$
exists. Namely, if this was the case, then Table~\ref{tab:8cycles} and Lemma~\ref{le:possHorbits} 
reveal that $v_a(\C) = 2 + 3 + v_a(\C_1) + v_a(\C_2)$ and
$v_g(\C) + v_z(\C) = 2 + 1 + v_g(\C_1) + v_z(\C_1) + v_g(\C_2) + v_z(\C_2)$. By Corollary~\ref{cor:summs}
and the fact that $v_a(\C) = v_g(\C) = v_z(\C)$, however, we have
$2(5 + v_a(\C_1) + v_a(\C_2)) = 3 + v_a(\C_1) + 3v_a(\C_2)$ and thus
$7 + v_a(\C_1) = v_a(\C_2)$. But since $v_a(\C_2) \leq 5$ (see Corollary~\ref{cor:summs}), this is impossible.
Thus no $8$-cycle of trace $a^5nan$ exists, as claimed.
Therefore, $v_a(\C) = 2 + v_a(\C_1) + v_a(\C_2)$ and 
$v_g(\C) + v_z(\C) = 2 + v_a(\C_1) + 3v_a(\C_2) $, so
$2 + v_a(\C_1) = v_a(\C_2)$. Since $v_a(\C_1)$ is even and $v_a(\C_2) \leq 5$, this leaves us with two
possible cases.

Case $v_a(\C_1) = 0$ and $v_a(\C_2) = 2$. Since no generic $8$-cycle and no $8$-cycle of $\C_2$ contains two
consecutive anchors, Lemma~\ref{le:quartic} implies that $X$ is half-arc-transitive.

Case $v_a(\C_1) = 2$ and $v_a(\C_2) = 4$. We show that this case actually cannot occur.
Note that, in view of Lemma~\ref{le:possHorbits}, we have 
$r^3 = \pm 1$ and no $8$-cycle corresponding to rows $22$-$29$ of Table~\ref{tab:8cycles} exists. 
Moreover, since $v_a(\C_1) = 2$, the set $\C_1$ consists precisely of one $H$-orbit of $8$-cycles which
corresponds to one of the rows $10$-$17$ of Table~\ref{tab:8cycles}. 
However, since $r^3 = \pm 1$, each of the corresponding conditions forces
some other condition of these rows to hold as well (just multiply by $r$ or $r^2$). This completes the proof.
\end{proof}

\begin{lemma}
\label{le:general}
Let $X = \X_e(m,n;r,t)$, where $m,n \geq 4$ are even integers and $r\in \ZZ_n^*$, 
$t \in \ZZ_n$ are such that  (\ref{eq:Xparamcond}) holds, and furthermore  $r^2 \neq \pm 1$.
If $X$ contains no $8$-cycle of trace $n^8$, $a^2n^6$, $anan^5$, $a^4n^4$, $a^2na^2n^3$ or $a^2n^2a^2n^2$, 
then $X$ is half-arc-transitive. In particular, if $m > 8$, then $X$ is half-arc-transitive.
\end{lemma}

\begin{proof}
Let $\C$ denote the set of all $8$-cycles of $X$. In view of Proposition~\ref{pro:halftrcond},
Proposition~\ref{pro:theisomorphisms} and Lemma~\ref{le:genequal} we can assume
that $v_a(\C) = v_g(\C)$ and $v_a(\C) \neq v_z(\C)$. Moreover, Proposition~\ref{pro:arc-tr} and
Proposition~\ref{pro:equalfreq} imply that if
$X$ is arc-transitive, then for any $\Aut X$-orbit of $8$-cycles $\C'$, we have $v_a(\C') = v_g(\C')$.
Let now $\C'$ denote the $\Aut X$-orbit of $8$-cycles of $X$ containing the generic $8$-cycles.
By Proposition~\ref{pro:equalfreq} every automorphism of $X$ maps a zigzag to a zigzag. Therefore, since
the generic $8$-cycles have precisely two zigzags which are antipodal, the same holds for all $8$-cycles
of $\C'$. Table~\ref{tab:8cycles} thus shows, 
that the only $H$-orbits of $8$-cycles that can be contained in $\C'$
are those of rows $1$, $10$, $17$, $18$ and $25$. 
Since $v_a(\C') = v_g(\C')$, at least one of 
the rows $18$ and $25$ corresponds to an $H$-orbit of $\C'$. 

Suppose first that the $H$-orbit of row $18$
is not in $\C'$ (and so the one of row $25$ is). Let $\C_1'$ denote the set of $8$-cycles of 
traces $a^2nan^2an$ and $a^3n^2an^2$ that are in $\C'$. We have 
$2 + v_a(\C_1') + 1 = v_a(\C') = v_g(\C') = 1 + v_g(\C_1') + 2$, that is $v_a(\C_1') = v_g(\C_1')$,
which implies (in view of the facts from Table~\ref{tab:8cycles}), that $\C_1'$ is empty. 
Therefore Lemma~\ref{le:quartic} implies that $X$ is half-arc-transitive. 

We can thus assume that $\C'$ contains the $H$-orbit corresponding to row $18$ of Table~\ref{tab:8cycles}
and so $r^3 = 1$. This implies that the $H$-orbit of row $25$ cannot exist for otherwise
$r-r^2 = 0$, that is $r=1$, which is impossible. In view of $v_a(\C') = v_g(\C')$, $\C'$ contains 
precisely one of the $H$-orbits of rows $10$ and $17$. However, in view of $r^3 = 1$, the condition of
row $10$ holds if and only if the condition of row $17$ holds. Therefore, the $H$-orbits of rows
$10$ and $17$ both exist and precisely one lies in $\C'$. Let $\C''$ denote the $\Aut X$-orbit of $8$-cycles
containing the other of the two $H$-orbits. As above every $8$-cycle of $\C''$ has precisely
two zigzags which are antipodal. But in view of the above remarks on such $8$-cycles, it is now clear that
$\C''$ consists of a single $H$-orbit of $8$-cycles of $X$, and so $2 = v_a(\C'') \neq v_g(\C'') = 1$.
Thus $X$ is half-arc-transitive, as claimed.
\end{proof}


\section{The case $m = 8$}

\indent
Throughout this section we let $m=8$ and $n,r,t,X,H,\rho,\sigma$ and $\tau$ be as in Section~\ref{sec:notation}.
Recall that this implies that the conditions (\ref{eq:Xparamcond}) are satisfied. In view of Proposition~\ref{pro:r^2}
we can assume $r^2 \neq \pm 1$.
In this section we show that $X$ is half-arc-transitive in this case. Recall (see Section~\ref{sec:8cycles})
that, apart from the $8$-cycles of Table~\ref{tab:8cycles}, the only possible $8$-cycles of $X$ are the
coiled ones, that is those of trace $n^8$. We now investigate the possible $H$-orbits 
of such $8$-cycles.
Note that we have
\begin{equation}
\label{eq:cond8}
	r^8 = 1\quad\quad \mathrm{and}\quad\quad 1+r+r^2 + \cdots + r^7 + 2t = 0.
\end{equation}
Observe also that every $8$-cycle of trace $n^8$ corresponds to a condition of the form
\begin{equation}
\label{eq:coiledcond}
	\delta_0 + \delta_1 r + \delta_2 r^2 + \cdots + \delta_7 r^7 + t = 0,\quad \mathrm{where}\
	\delta_i \in \{0,1\}\ \mathrm{for}\ \mathrm{all}\ i \in \{0,1,\ldots , 7\}.
\end{equation}
Moreover, since we are only interested
in the $H$-orbits of such $8$-cycles, we can assume that at most $4$ of the numbers $\delta_i$ are nonzero.
The next lemma gives all the possible $H$-orbits of coiled $8$-cycles of $X$. 

\begin{lemma}
\label{le:coiled8}
Let $X = \X_e(m,n;r,t)$, where $m = 8$, $n \geq 4$ is even and $r\in \ZZ_n^*$, 
$t \in \ZZ_n$ are such that  (\ref{eq:Xparamcond}) holds, and furthermore  $r^2 \neq \pm 1$.
Then the only possible $H$-orbits of coiled $8$-cycles of $X$ are those listed in Table~\ref{tab:coiled8}.
Moreover, the following hold:
\begin{itemize}
	\item[(i)] The $H$-orbit of row $3$ exists if and only if the $H$-orbit of row $6$ exists.
	\item[(ii)] If the $H$-orbit of row~2 exists, then no other $H$-orbit of coiled $8$-cycles exists.
	\item[(iii)] At most one $H$-orbit of rows $3$, $4$ and $5$ can exist.
	\item[(iv)] If the $H$-orbit of row $4$ exists, then either all or none of the $H$-orbits of rows $1$, $7$ and $8$ exists.
\end{itemize}
\end{lemma}

\begin{table}[!hbt]
\begin{footnotesize}
\begin{center}
\begin{tabular}{@{}c|l|c|c|c@{}}
	Row	&  A representative 	& Code	& Condition	& Orbit length \\
	\hline  & & & & \\
   $1$ & $u_0^0u_1^0u_2^0u_3^0u_4^0u_5^0u_6^0u_7^0$ & $g^8$ & $t = 0$ & $2n$ \\    
   $2$ & $u_0^0u_1^0u_2^0u_3^0u_4^0u_5^0u_6^{r^5}u_7^{r^5}$ & $g^4z^4$ & $1+r^2+t=0$ & $16n$ \\
   $3$ & $u_0^0u_1^0u_2^0u_3^0u_4^0u_5^{r^4}u_6^{r^4+r^5}u_7^{r^4+r^5+r^6}$ & $g^3zg^3z$ & $1+r+r^2+r^3+t=0$ & $8n$ \\
   $4$ & $u_0^0u_1^0u_2^0u_3^0u_4^{r^3}u_5^{r^3}u_6^{r^3}u_7^{r^3}$ & $g^2z^2g^2z^2$ & $1+r^4+t=0$ & $8n$ \\
   $5$ & $u_0^0u_1^0u_2^0u_3^0u_4^{r^3}u_5^{r^3+r^4}u_6^{r^3+r^4}u_7^{r^3+r^4+r^6}$ & $g^2zgz^2gz$ & $1+r+r^3+r^4+t=0$ & $16n$ \\
   $6$ & $u_0^0u_1^0u_2^0u_3^{r^2}u_4^{r^2}u_5^{r^2+r^4}u_6^{r^2+r^4+r^5}u_7^{r^2+r^4+r^5}$ & $gz^3gz^3$ & $1+r^2+r^3+r^5+t=0$ & $8n$ \\
   $7$ & $u_0^0u_1^0u_2^0u_3^{r^2}u_4^{r^2+r^3}u_5^{r^2+r^3}u_6^{r^2+r^3}u_7^{r^2+r^3+r^6}$ & $gzgzgzgz$ & $1+r+r^4+r^5+t=0$ & $4n$ \\
   $8$ & $u_0^0u_1^0u_2^ru_3^{r}u_4^{r+r^3}u_5^{r+r^3}u_6^{r+r^3+r^5}u_7^{r+r^3+r^5}$ & $z^8$ & $1+r^2+r^4+r^6+t=0$ & $2n$
\end{tabular}
\caption{Possible $H$-orbits of coiled $8$-cycles of $X$ when $m=8$.}
\label{tab:coiled8}
\end{center}
\end{footnotesize}
\end{table}

\begin{proof}
As already noted each $H$-orbit of coiled $8$-cycles corresponds to a condition as in $(\ref{eq:coiledcond})$ 
and we can assume that at most four of the numbers
$\delta_i$ are nonzero. We distinguish five cases depending on the number $s$ of nonzero 
multipliers $\delta_i$.\smallskip

Case $s = 0$. The condition is then $t=0$ and the $H$-orbit in question is clearly the one 
corresponding to row~1 of Table~\ref{tab:coiled8}.\smallskip

Case $s = 1$. We show that this is not possible. Namely, in this case the condition is $r^i + t = 0$ for some 
$i \in \{0,1,\ldots , 7\}$. Thus, in view of $t(r-1) = 0$, we have $r^i(r-1) = 0$,
and so $r-1 = 0$, a contradiction.\smallskip

Case $s = 2$. The condition is then of the form $1+r^i + t = 0$ for some
$i \in \{1,2,\ldots , 7\}$. Therefore $r^i + r^{2i} + t=0$ and thus
$r^{2i}-1 = 0$. In view of $r^2 \neq \pm 1$, this implies $2i = 0$ or $2i = 4$. We can thus assume
$i = 4$ or $i = 2$. In the first case $1+r^4+t=0$ and the $H$-orbit in question corresponds to
row $4$ of Table~\ref{tab:coiled8}, and in the second case $1+r^2+t=0$ and the $H$-orbit
corresponds to row $2$ of Table~\ref{tab:coiled8}.\smallskip

Case $s = 3$. We show that this also is not possible. Namely, any such condition is of the form
$1+r^i + r^j + t=0$, where $1 \leq i < j \leq 7$. It is easy to see that we can multiply
this equation by an appropriate power $r^k$, so that $0,i,j,k,i+k$ and $j+k$ are all distinct modulo $8$.
Then $1+r^i+r^j+r^k+r^{i+k}+r^{j+k}+2t = 0$ and (\ref{eq:cond8}) imply, that $1+r^l=0$ 
for some $l \in \{1,2,\ldots , 7\}$. In view of $r^8 = 1$ and $r^2 \neq \pm 1$, we have
$l = 4$, and so $r^4 = -1$. This implies that none of $i,j$ or $j-i$ equals $4$, since otherwise
a condition of the form $r^{k'}+t=0$ is obtained, which is impossible. With no loss of generality we
can assume that $i$ is minimal among $i, j-i, 8-j$. Then the pair $(i,j)$ is one of
$(1,2)$, $(1,3)$, $(1,6)$, $(1,7)$ or $(2,5)$. Note that $(1,7)$ is equivalent to $(1,2)$ in the sense that
multiplying $1+r+r^7+t = 0$ by $r$ we get $1+r+r^2+t=0$. 
If $(i,j) = (1,2)$, then $1+r+r^2+t=0$, so
$r^5+r^6+r^7+t=0$ and thus (\ref{eq:cond8}) forces $r^3+r^4=0$, giving $r=-1$, a contradiction. A similar
contradiction is obtained if $(i,j) = (2,5)$. If $(i,j) = (1,3)$, however, then
$1+r+r^3+t=0$ and so since $r^4 = -1$, multiplication by $r$ gives $-1+r+r^2+t=0$. Subtracting we have $r^3-r^2+2=0$. 
Multiplying by $r$ we get $-1-r^3+2r=0$, and so $0 = r^3-r^2+2+(-1-r^3+2r) = -r^2+2r+1$. Adding this to 
$-1+r+r^2+t=0$ finally forces $3r+t=0$. Thus $3+t=0$ and so $t = -3$. But $r^4=-1$ and (\ref{eq:cond8})
imply $2t=0$, and so $6=0$, contradicting (\ref{eq:n}). A similar contradiction is obtained in case
$(i,j) = (1,6)$.\smallskip

Case $s = 4$. Then the condition is of the form $1+r^i+r^j+r^k+t=0$, where $1 \leq i < j < k \leq 7$. 
Note that the minimum of
numbers $i,j-i,k-j,8-k$ is either $1$ or $2$. In the latter case the condition is clearly
$1+r^2+r^4+r^6+t=0$ and the $H$-orbit in question corresponds to row~8 of Table~\ref{tab:coiled8}. 
In the former case
it is easy to see that (using $(\ref{eq:cond8})$) we can assume that the triple $(i,j,k)$ is one of 
$(1,2,3)$, $(1,2,4)$, $(1,2,5)$, $(1,3,5)$, $(1,3,6)$
or $(1,4,5)$. The triples $(1,2,3)$, $(1,2,5)$, $(1,3,6)$ and $(1,4,5)$ give rise to $H$-orbits of
rows $3$, $5$, $6$ and $7$, respectively. We show that $(1,2,4)$ and $(1,3,5)$ are impossible.
The triple $(1,2,4)$ gives rise to the condition $1+r+r^2+r^4+t=0$. Multiplying by $r^5$ we get
$r^5+r^6+r^7+r+t=0$, so adding the two equations we get $1+2r+r^2+r^4+r^5+r^6+r^7+2t=0$. 
By (\ref{eq:cond8}), however, we
have $r-r^3=0$,
that is $r^2 = 1$, a contradiction. A similar contradiction is obtained by the triple $(1,3,5)$.\smallskip

This proves that the only possible $H$-orbits of coiled $8$-cycles of $X$ are those listed in Table~\ref{tab:coiled8}.
Representatives, codes and orbit lengths are now easily obtained.
As for the second part, the four claims, we proceed as follows. Observe first that if $r^4 = 1$, then
each of the conditions of rows $3$ and $6$ forces the other to hold as well. Now if
$1+r+r^2+r^3+t=0$, then $r+r^2+r^3+r^4+t=0$, so subtracting the two equations we get $r^4-1=0$, as desired. If however
$1+r^2+r^3+r^5+t=0$, then $r+r^3+r^4+r^6+t=0$, and so (\ref{eq:cond8}) implies $r^7-r^3=0$, that is
$r^4=1$. This proves the first claim. 

Suppose now that the condition of row~2 holds, that is 
$1+r^2+t=0$. Then $r+r^3+t=0$, and so (\ref{eq:cond8}) gives $r^4+r^5+r^6+r^7=0$, that is
$1+r+r^2+r^3 = 0$. This forces $2t=0$, so since $1+r^2 \neq 0$, we have $t = \frac{n}{2}$. 
Moreover, since $1+r^2$ is even (recall that $r$ is odd), $t$ is also even, so
$n \equiv 0 \pmod 4$. Thus $r^4 = (r^2+1-1)^2 = (\frac{n}{2})^2 - n + 1 = 1$. Clearly,
the conditions of rows $1$ and $4$ are impossible in view of (\ref{eq:n}). Moreover, the
condition of row~3 (and thus by the preceding paragraph also of row~6) also cannot hold. 
If $1+r+r^3+r^4+t=0$, then
$0 = 1+r\frac{n}{2} + 1 + \frac{n}{2} =  2$, a contradiction. If $1+r+r^4+r^5+t=0$,
then $r+r^2+r^5+r^6+t=0$, and so $0=r^6-r^4+r^2-1=(r^4+1)(r^2-1) = 2(\frac{n}{2}-2) = -4$, a contradiction.
Finally, if $1+r^2+r^4+r^6+t=0$, then $r^4\frac{n}{2} = \frac{n}{2} = 0$, a contradiction. 
Thus the only $H$-orbit of coiled $8$-cycles is that of row $2$, as claimed.

Simultaneous existence of $H$-orbits of rows $3$ and $5$ or $4$ and $5$ would
contradict $r^2 \neq \pm1$. As for $3$ and $4$, it was shown above that the condition of 
row $3$ implies $r^4 = 1$, so if $1+r^4+t=0$, we get $t=-2$, which, in view of $t(r-1)=0$, contradicts
(\ref{eq:2r}). This proves the third claim.

The last claim is straightforward.
\end{proof}

\begin{lemma}
\label{le:m=8equal}
Let $X = \X_e(m,n;r,t)$, where $m = 8$, $n \geq 4$ is an even integer and $r\in \ZZ_n^*$, 
$t \in \ZZ_n$ are such that  (\ref{eq:Xparamcond}) holds, and furthermore  $r^2 \neq \pm 1$.
If $v_a(8) = v_g(8) = v_z(8)$, then $X$ is half-arc-transitive.
\end{lemma}

\begin{proof}
Note that the assumptions imply that all the possible $H$-orbits of $8$-cycles of $X$ are those listed in
Table~\ref{tab:8cycles} and Table~\ref{tab:coiled8}.
Moreover, since $n \geq 14$ is even and since there is no element $r$ of
$\ZZ_{14}^*$, such that $r^8 = 1$ and $r^2 \neq 1$, we can assume $n \geq 16$.
Let $\C$ denote the set of all $8$-cycles of $X$, 
let $\C_1$ denote the set of $8$-cycles of $X$ of traces $a^2nan^2an$ and $a^3n^2an^2$, let $\C_2$
denote the set of $8$-cycles of $X$ of trace $an^3an^3$ and let $\C_c$ denote the set of 
coiled $8$-cycles of $X$. In view of Lemma~\ref{le:genequal} we can assume that $\C_c$ is nonempty.

We claim that no $8$-cycle of trace $a^5nan$ exists in $X$. If this is not the case, then, in view of 
Proposition~\ref{pro:theisomorphisms}, we can assume $r=3$. Since $n \geq 16$, the only $H$-orbits that
could be contained in $\C_1$ are those corresponding to rows $6$, $14$, $15$ and $16$ of Table~\ref{tab:8cycles}.
Note however, that in case of row~14 we have $22=0$, but then $r^8 = 5 \neq 1$, which is impossible.
Thus either $n = 16$, in which case $\C_1$ consists of the $H$-orbits corresponding to rows $6$ and $16$ of
Table~\ref{tab:8cycles}, or $n = 20$, in which case $\C_1$ consists of the $H$-orbit corresponding to row $15$ of
Table~\ref{tab:8cycles}, or $\C_1$ is empty.
Suppose first that $n = 16$. It is easy to check that the conditions of rows $25$ and $28$ of Table~\ref{tab:8cycles}
hold, and so Lemma~\ref{le:possHorbits} implies that $\C_2$ consists precisely of the two corresponding $H$-orbits. 
Since $v_a(\C_c)=0$, we thus have $v_a(\C) = 2 + 3 + 6 + 2 = 13$. 
Moreover, $v_g(\C) = 1 + 0 + 4 + 3 + v_g(\C_c)$ and
$v_z(\C) = 1 + 1 + 2 + 3 + v_z(\C_c)$, so, in view of $v_a(\C) = v_g(\C) = v_z(\C)$, we have $v_g(\C_c) = 5$ and
$v_z(\C_c) = 6$. 
Since $v_g(\C_c) + v_z(\C_c) = 11$ is odd, the data of
Table~\ref{tab:coiled8} reveals that precisely one $H$-orbit corresponding to rows $1$ and $8$
of Table~\ref{tab:coiled8} exists.
Note that $r^4 = 1$ and $1+r+r^2+r^3 = 8$. Thus equation
(\ref{eq:cond8}) implies $2t = 0$. 
It follows that $1+r^2+r^4+r^6+t = 4+t$ cannot be zero and so the $H$-orbit of row~1 of Table~\ref{tab:coiled8}
exists, that is $t = 0$, and so $X = \X_e(8,16;3,0)$.
It is now easy to check that $C_c$ consists precisely of the
$H$-orbits corresponding to rows $1$ and $5$ of Table~\ref{tab:coiled8}, and so $v_z(\C_c) = 4 \neq 6$, a contradiction.
Thus $n$ cannot be $16$.
Suppose now that $n = 20$. Then $\C_1$ consists solely of the $H$-orbit corresponding to 
row $15$ of Table~\ref{tab:8cycles}. Moreover, as the conditions of rows $22$ and $27$ of Table~\ref{tab:8cycles}
both hold, Lemma~\ref{le:possHorbits} implies that $\C_2$ consists precisely of the two corresponding
$H$-orbits. We thus have
$v_a(\C) = 2 + 3 + 2 + 2 = 9$, $v_g(\C) = 1 + 0 + 1 + 3 + v_g(\C_c)$ and $v_z(\C) = 1 + 1 + 1 + 3 + v_z(\C_c)$.
Therefore, $v_g(\C_c) = 4$ and $v_z(\C_c) = 3$. As above precisely one $H$-orbit corresponding to
rows $1$ and $8$
of Table~\ref{tab:coiled8} exists. But this is impossible since $1 + r^2 + r^4 + r^6 = 0$. We are thus left with
the possibility that $\C_1$ is empty. Using Corollary~\ref{cor:summs}, we then get 
$v_a(\C) = 2 + 3 + v_a(\C_2)$ and 
$v_g(\C) + v_z(\C) = 2 + 1 + 3v_a(\C_2) + v_g(\C_c) + v_z(\C_c)$, and so 
$7 = v_a(\C_2) + v_g(\C_c) + v_z(\C_c)$. Since $r^3 = \pm 1$ clearly cannot hold, 
Lemma~\ref{le:possHorbits} implies that
$v_a(\C_2) \leq 2$, and so $5 \leq v_g(\C_c) + v_z(\C_c) \leq 7$. Together with Lemma~\ref{le:coiled8}, this implies
that no $H$-orbit of rows $2$, $3$, $5$ and $6$ of Table~\ref{tab:coiled8} exists. Moreover, the $H$-orbit 
corresponding to row $4$ of Table~\ref{tab:coiled8} also cannot exist. Namely, in view of
$v_g(\C_c) + v_z(\C_c) \geq 5$, at least one $H$-orbit out of those corresponding to rows $1$, $7$ and $8$ would 
also exist. By Lemma~\ref{le:coiled8}, they would all exist, forcing 
$v_g(\C_c) + v_z(\C_c) = 8$, a contradiction.
Therefore, $\C_c$ consist solely of some of the $H$-orbits corresponding to 
rows $1$, $7$ and $8$ of Table~\ref{tab:coiled8}.
But then $v_g(\C_c) + v_z(\C_c) \leq 4$, which is impossible. 
Thus no $8$-cycle of trace $a^5nan$ exists, as claimed.

In view of Lemma~\ref{le:quartic}, we can assume that $\C_1$ is nonempty. 
But then Lemma~\ref{le:C1orb} implies that $\C_1$ consists of precisely one $H$-orbit of $8$-cycles. 
Namely, since $r^8 = 1$, the order of $r$ in $\ZZ_n^*$ cannot be $3$, $5$ or $6$.
We now distinguish two possible cases.

Suppose first that $v_a(\C_1) = 4$, that is $\C_1$ is an $H$-orbit of $8$-cycles of trace $a^2nan^2an$.
Since $2v_a(\C) = v_g(\C) + v_z(\C)$, Corollary~\ref{cor:summs} implies
$$ 2(2+4+v_a(\C_2)) = 2 + 4 + 3v_a(\C_2) + v_g(\C_c) + v_z(\C_c),\ \mathrm{so}\ v_a(\C_2) + v_g(\C_c) + v_z(\C_c) = 6. $$
In view of Lemma~\ref{le:coiled8}, this implies that none of the rows $2,3,5$ or $6$ of
Table~\ref{tab:coiled8} corresponds to an $H$-orbit of $8$-cycles of $X$. Furthermore, 
row $4$ of Table~\ref{tab:coiled8}
cannot correspond to an $H$-orbit of $8$-cycles of $X$. Namely, if this was the case, then Lemma~\ref{le:coiled8} and
$v_g(\C_c) + v_z(\C_c) \leq 6$ imply, that $\C_c$ consists precisely of that $H$-orbit and thus 
$v_a(\C_2) = 2$. But then Lemma~\ref{le:possHorbits} implies $r^2 = \frac{n}{2} \pm 1$, so
$r^4 = 1$, which together with $1+r^4+t=0$ forces $t=-2$. But in view of $(\ref{eq:Xparamcond})$, this
contradicts (\ref{eq:2r}).
Clearly $r^3 \neq \pm 1$, so we
have $v_a(\C_2) \leq 2$, and thus $4 \leq v_g(\C_c) + v_z(\C_c) \leq 6$.
But then $v_g(\C_c) + v_z(\C_c) = 4$, and so $\C_c$ consists precisely of the
$H$-orbits corresponding to rows $1$, $7$ and $8$ of Table~\ref{tab:coiled8}. Thus $v_a(\C_2) = 2$,
so as above $r^4 = 1$. Together with $t=0$ and $1+r+r^4+r^5+t=0$, this gives $2(1+r)=0$, which
contradicts $(\ref{eq:2r})$.

Suppose now that $v_a(\C_1) = 2$. Similarly as in the previous paragraph we have
\begin{equation}
\label{eq:8aux}
	2(2+2+v_a(\C_2)) = 2 + 2 + 3v_a(\C_2) + v_g(\C_c) + v_z(\C_c),\ \mathrm{so}\ v_a(\C_2) + v_g(\C_c) + v_z(\C_c) = 4. 
\end{equation}
Thus none of the rows $2,3,5$ or $6$ of Table~\ref{tab:coiled8} corresponds to an $H$-orbit of $8$-cycles of $X$.
In view of Proposition~\ref{pro:theisomorphisms}, we can assume that precisely
one of the $H$-orbits of rows $12$ and $13$ of Table~\ref{tab:8cycles} exists. 

Suppose first that
$\C_1$ is the $H$-orbit corresponding to row $13$ of Table~\ref{tab:8cycles}, that is $2+r+r^2=0$.
Then $r^4 = (-2-r)^2 = 4+4r+r^2=2+3r$ and thus $r^8 = (2+3r)^2 = 4 + 12r + 9r^2 = -14 + 3r$.
Therefore, $r^8=1$ implies $3(r-5)=0$, and so $0 = 6 + 3r + 3r^2 = 96$,
that is $96 \equiv 0 \pmod n$.
We claim that $3$ does not divide $n$. Namely, if it does, then in view of $r \in \ZZ_n^*$, we
either have $r \equiv 1 \pmod 3$ or $r \equiv 2 \pmod 3$. But in the former case
$2+r+r^2 \equiv 1 \pmod 3$ and in the latter case $2+r+r^2 \equiv 2 \pmod 3$, contradicting $2+r+r^2=0$.
Thus $3$ does not divide $n$, and so $r = 5$. Then $2+r+r^2=0$ implies $32 \equiv 0 \pmod n$, so
in view of the fact that no $8$-cycle of trace $a^5nan$ exists, $n = 32$.
This implies $1+r+\cdots + r^7 = 24$, so
$2t = 8$. But $(\ref{eq:Xparamcond})$ implies $0 = t\cdot 4 = 2t\cdot 2 = 16$, a contradiction.

We can now finally assume that 
$\C_1$ is the $H$-orbit corresponding to row $12$ of Table~\ref{tab:8cycles}, 
that is $2+r-r^2=0$. We thus have $r^4=6+5r$, so
$r^8 = 86 + 85r$. Together with $r^8 = 1$, this implies $5\cdot 17 (r+1) = 0$. Therefore,
at least one of $5$ and $17$ divides $n$.
Observe first that $v_a(\C_2) = 2$ is impossible, since then
Lemma~\ref{le:possHorbits} implies that one of the conditions in rows $25$ or $27$ of Table~\ref{tab:8cycles}
holds. But the one of row $25$ gives $r^3 = -1$ and the one of row $27$ gives $1-3r=0$, which are
both impossible. Suppose now that $v_a(\C_2) = 1$. Then $(\ref{eq:8aux})$ implies 
$v_g(\C_c) + v_z(\C_c) = 3$. Moreover,
$5 = v_a(\C) = v_g(\C) = 1 + 1 + v_g(\C_2) + v_g(\C_c)$ and $5 = v_z(\C) = 1+1+v_z(\C_2)+v_z(\C_c)$. 
Thus, in view of the facts from Table~\ref{tab:coiled8}, 
precisely one of $v_g(\C_2)$ and $v_z(\C_2)$ is $1$ and the other is $2$. 
Therefore, $\C_2$ is the $H$-orbit corresponding
to one of the rows $23$, $24$, $25$, $26$, $28$ or $29$ of Table~\ref{tab:8cycles}. 
However, all but the condition of row $29$ are impossible. Namely,  
row $23$ gives $0 = 1+r+r^2-r^3 = 1 + r - 2r = 1-r$, row $24$ gives $1+3r=0$, row $25$ gives $r^3=-1$, 
row $26$ gives
$1+r+2r^2=0$ and row $28$ gives $1+r=0$. These are all impossible (recall that $1+r+2r^2=0$ cannot hold
since $\C_1$ consists solely of the $H$-orbit corresponding to row~12 of Table~\ref{tab:8cycles}). 
Thus the condition of row $29$ holds, that is $1-r-r^2-r^3=0$.
Since $r^3 = 2r+r^2=2+3r$, we thus have $-3-5r=0$. Therefore, $0 = 10+5r-5r^2=7+3r$. 
Hence $4-2r=0$ and thus also $11+r=0$. But then $n = 26$ which
contradicts the fact that at least one of $5$ and $17$ divides $n$. 
We are left with the possibility that $v_a(\C_2) = 0$.
Thus $v_g(\C_c) + v_z(\C_c) = 4$ and $v_g(\C_c) = v_z(\C_c) = 2$. By Lemma~\ref{le:coiled8},
either $\C_c$ consists of the $H$-orbit corresponding to row $4$ of Table~\ref{tab:coiled8}, or
it consists of the $H$-orbits corresponding to rows $1$, $7$ and $8$ of that table. We consider each case
separately.

If $\C_c$ consists of the $H$-orbit corresponding to row $4$ of Table~\ref{tab:coiled8},
then $1+r^4+t=0$ and thus also $r+r^5+t=0$. In view of (\ref{eq:cond8}), we have $r^2+r^3+r^6+r^7=0$.
Thus $1+r+r^4+r^5=0$, which together with $1+r^4+t+r+r^5+t=0$ gives $2t=0$. 
Then $t = \frac{n}{2}$, since $t \neq 0$ in view of the fact that the
$H$-orbit of row $1$ from Table~\ref{tab:coiled8} does not exist. Therefore,
$2+2r^4=0$, so since $r^4 = 6+5r$, we have $14+10r=0$. 
But $0 = 10(2+r-r^2) = -8+4r$,
which together with $14+10r=0$ gives $30 + 2r=0$. Thus $60+4r=0$, and so $68 \equiv 0 \pmod n$.
Since $5 \cdot 17(r+1) = 0$, we thus have $17(r+1)=0$, and so $r+49=0$. If $n = 34$, then
$r = 19$, and so since $t = \frac{n}{2} = 17$, we have $1+r^4+t = 1 + 33+ 17 = 17 \neq 0$, which is impossible.
Thus $n = 68$ and $r = 19$. Moreover, $t = \frac{n}{2} = 34$, and so $X = \X_e(8,68;19,34)$.
We show that $X$ is half-arc-transitive in this case. Suppose on the
contrary that it is arc-transitive. Since no
generic and no coiled $8$-cycle has two consecutive anchors, Lemma~\ref{le:quartic} implies
that the $8$-cycles of the $H$-orbit corresponding to row~12 of Table~\ref{tab:8cycles}
are in the same $\Aut X$-orbit as the generic $8$-cycles. However,
each of the $3$-paths of the $8$-cycle 
$u_0^0u_1^0u_2^0u_3^{r^2}u_2^{r^2}u_1^{-r+r^2}u_0^{-1-r+r^2}u_1^{-1-r+r^2}$, denoted by $C_1$,
except $u_1^{-r+r^2}u_0^{-1-r+r^2}u_1^{-1-r+r^2}u_0^0$ and $u_0^{-1-r+r^2}u_1^{-1-r+r^2}u_0^0u_1^0$ lie
on two $8$-cycles of $X$. Consider on the other hand the generic $8$-cycle $u_0^0u_1^0u_2^0u_1^{-r}u_0^{-r}u_1^{1-r}u_2^1u_1^1$
and denote it by $C_2$.
It turns out that every other $3$-path of $C_2$ lies on two $8$-cycles of $X$ and every other $3$-path of $C_2$ on
one $8$-cycle of $X$. Therefore, Proposition~\ref{pro:equalfreq} implies 
that no automorphism of $X$ can map $C_1$ to $C_2$ and so
they are not in the same $\Aut X$-orbit, a contradiction. Thus $X$ is half-arc-transitive, as claimed.

If on the other hand $\C_c$ consists of the $H$-orbits corresponding to rows $1$, $7$ and $8$ of Table~\ref{tab:coiled8},
then the argument is as follows. Since $r^4 = 6+5r$, we have $r^5 = 10 + 11r$ and $r^6 = 22 + 21r$.
Thus, since $t = 0$, the condition of row $7$ of Table~\ref{tab:coiled8} implies
$17 + 17r = 0$. Moreover, the condition of row $8$ implies
$31+27r=0$. Therefore, $3+7r=0$ and consequently $11+3r=0$. This gives $-19+r=0$,
and so $68\equiv 0 \pmod n$. If $n = 34$, then $r^4 + 1 = 0$ which contradicts the fact that row $4$
does not correspond to an $H$-orbit of $8$-cycles of $X$. Therefore, $X = \X_e(8,68;19,0)$.
As in the previous paragraph one can check that on a generic $8$-cycle every other $3$-path lies on two
$8$-cycles of $X$ and every
other $3$-path lies on one $8$-cycle of $X$, whereas on an $8$-cycle of the $H$-orbit corresponding
to row~12 of Table~\ref{tab:8cycles} all but two $3$-paths lie on two $8$-cycles of $X$, and so $X$ is
half-arc-transitive.
\end{proof}

\begin{lemma}
\label{le:m=8}
Let $X = \X_e(m,n;r,t)$, where $m = 8$, $n \geq 4$ is an even integer and $r\in \ZZ_n^*$, 
$t \in \ZZ_n$ are such that  (\ref{eq:Xparamcond}) holds, and furthermore  $r^2 \neq \pm 1$.
Then $X$ is half-arc-transitive.
\end{lemma}

\begin{proof}
As in the proof of Lemma~\ref{le:m=8equal}, we can assume $n \geq 16$.
Let $\C$ denote the set of all $8$-cycles of $X$. In view of Lemma~\ref{le:general} 
we can assume that $\C$ contains at least one coiled $8$-cycle.
As in the proof of Lemma~\ref{le:general} we can assume
that $v_a(\C) = v_g(\C)$ and $v_a(\C) \neq v_z(\C)$. Moreover, if 
$X$ is arc-transitive, then for any $\Aut X$-orbit of $8$-cycles $\C'$ we have $v_a(\C') = v_g(\C')$.
Let now $\C'$ denote the $\Aut X$-orbit of $8$-cycles of $X$ containing the generic $8$-cycles.
By Proposition~\ref{pro:equalfreq}, every automorphism of $X$ maps a zigzag to a zigzag. Thus since
the generic $8$-cycles have precisely two zigzags which are antipodal, the same holds for all $8$-cycles
of $\C'$.

Suppose first that $\C'$ contains no coiled $8$-cycles. 
Since $r^3 \neq \pm 1$, the arguments of the proof of Lemma~\ref{le:general} show, that
$\C'$ contains no $8$-cycle with two consecutive anchors, and so Lemma~\ref{le:quartic} implies
that $X$ is half-arc-transitive.

Suppose now that $\C'$ contains some coiled $8$-cycles. Since they also have precisely two zigzags
which are antipodal, Table~\ref{tab:coiled8} shows, that the coiled $8$-cycles of $\C'$ are
precisely the ones of the $H$-orbit corresponding to row $3$ of Table~\ref{tab:coiled8}.
Let $\C_1'$ denote the set of $8$-cycles of $\C'$ of traces $a^2nan^2an$ and $a^3n^2an^2$
and let $\C_2'$ denote the set of $8$-cycles of $\C'$ of trace $an^3an^3$. The
only $H$-orbits of $8$-cycles that can be contained in $\C_1'$ are those corresponding to rows $10$ and $17$ 
of Table~\ref{tab:8cycles}. By Lemma~\ref{le:quartic} we can assume that $\C_1'$ is nonempty. 
This implies that 
precisely one of the two $H$-orbits lies in $\C_1'$, for otherwise
$2-r-r^2=0$ and $1+r-2r^2=0$, which implies
$1 - 2r + r^2=0$, that is $(r-1)^2=0$. But since also $3 - 3r=0$, we get
$r^3 = (r-1+1)^3=1$, a contradiction.
Therefore, $\C_1'$ consists of precisely one $H$-orbit of $8$-cycles, and so $v_a(\C_1') = 2$ and
$v_g(\C_1') = 1$. Since $v_a(\C') = v_g(\C')$, the equations $v_a(\C') = 2 + 2 + v_a(\C_2')$ and 
$v_g(\C') = 1 + 1 + v_g(\C_2') + 3$, imply that $v_a(\C_2') - v_g(\C_2') = 1$, which is impossible,
since $\C_2'$ is either empty or consists precisely of the $H$-orbit corresponding 
to row $25$ of Table~\ref{tab:8cycles}. Thus $X$ is half-arc-transitive.
\end{proof}


\section{The case $m = 6$}

\indent
Throughout this section we let $m=6$ and we let $n,r,t,X,H,\rho,\sigma$ and $\tau$ be as in Section~\ref{sec:notation}.
Recall that this implies that the conditions (\ref{eq:Xparamcond}) are satisfied. In view of Proposition~\ref{pro:r^2}
we can assume $r^2 \neq \pm 1$. Recall also that we can assume 
(\ref{eq:2r}) and (\ref{eq:n}). In view of the results of
Section~\ref{sec:8cycles}, the only possible $H$-orbits of $8$-cycles of $X$ are those listed in 
Table~\ref{tab:8cycles}
and the ones with trace $a^2n^6$. The following observations give us some information about the
$8$-cycles of trace $a^2n^6$.

\begin{lemma}
\label{le:a2n6}
Let $X = \X_e(m,n;r,t)$, where $m = 6$, $n \geq 4$ is an even integer and $r\in \ZZ_n^*$, 
$t \in \ZZ_n$ are such that  (\ref{eq:Xparamcond}) holds, and furthermore  $r^2 \neq \pm 1$.
Then the length of each $H$-orbit of $8$-cycles of trace $a^2n^6$ is $2mn$ and each $8$-cycle of trace $a^2n^6$
has an even number of zigzags and an even number of glides. 
\end{lemma}

\begin{proof}
Note that for each $8$-cycle $C$ of trace $a^2n^6$ there is a unique $i \in \ZZ_6$, so that
$C$ contains two vertices from each of the attachment sets $X_i$, $X_{i+1}$ and one vertex from 
each of the other four attachment sets. Denote the two vertices of $X_i$ by $u$ and $v$ and the two vertices
of $X_{i+1}$ by $w$ and $x$. Out of the latter two vertices precisely one is adjacent to both $u$ and $v$.
With no loss of generality assume it is $w$. Suppose now 
there exists a nontrivial automorphism $\varphi$ of $H$ fixing $C$ setwise.
Since the sets $X_j$, $j \in \ZZ_6$, are blocks of imprimitivity for $H$ and since
$H$ acts half-arc-transitively on $X$, we have $\{u,v\}\varphi = \{u,v\}$. Thus $w$ is left fixed by $\varphi$.
But $H_w$ is of order $2$, and so $\varphi$ is the unique nontrivial element of $H_w$, which implies that
$x$ is not fixed by $\varphi$. However, since $\{w,x\}\varphi = \{w,x\}$, we have $x\varphi = x$, 
a contradiction. Thus no element of $H$ fixes $C$ setwise
and thus the $H$-orbit of $C$ has length $2mn$, as desired. The second claim of the lemma 
follows from Proposition~\ref{pro:nonanchors}.
\end{proof}

\begin{lemma}
\label{le:m=6equal}
Let $X = \X_e(m,n;r,t)$, where $m = 6$, $n \geq 4$ is an even integer and $r\in \ZZ_n^*$, 
$t \in \ZZ_n$ are such that  (\ref{eq:Xparamcond}) holds, and furthermore  $r^2 \neq \pm 1$. 
If $v_a(8) = v_g(8) = v_z(8)$, then $X$ is half-arc-transitive.
\end{lemma}

\begin{proof}
Let $\C$ denote the set of all $8$-cycles of $X$, let $\C_1$ denote the set of $8$-cycles 
of traces $a^2nan^2an$ and $a^3n^2an^2$, let $\C_2$
denote the set of $8$-cycles of trace $an^3an^3$ and let $\C_3$ denote the set of 
$8$-cycles of trace $a^2n^6$. In view of Lemma~\ref{le:genequal} we can assume that $\C_3$ is nonempty.

We claim that no $8$-cycle of trace $a^5nan$ exists in $X$. Suppose this is not the case. In view of
Proposition~\ref{pro:theisomorphisms} we can assume $r=3$. 
In view of Corollary~\ref{cor:summs} we have 
$v_a(\C) = 2 + 3 + v_a(\C_1) + v_a(\C_2) + v_a(\C_3)$
and $v_g(\C) + v_z(\C) = 2 + 1 + v_a(\C_1) + 3v_a(\C_2) + 3v_a(\C_3)$. Thus $v_a(\C) = v_g(\C) = v_z(\C)$
implies $7 + v_a(\C_1) = v_a(\C_2) + v_a(\C_3)$. Since $v_a(\C_1)$ and $v_a(\C_3)$ are even (see
Corollary~\ref{cor:summs} and Lemma~\ref{le:a2n6}), $v_a(\C_2)$ is odd, which, in view of 
Lemma~\ref{le:possHorbits}, implies that precisely one of the rows
$22$-$29$ of Table~\ref{tab:8cycles} 
corresponds to an $H$-orbit of $8$-cycles of $X$. It is easy to check that
each of the corresponding conditions, except for that of row $23$, contradicts the fact that $n \geq 14$ and $r^6 = 1$.
We can thus assume that the condition of row $23$ holds, that is $1+r+r^2-r^3=0$, and so $n = 14$. 
This implies $r^3 = -1$, and so $v_a(\C_2) = 5$. Moreover, out of the conditions of rows $10$-$17$ 
of Table~\ref{tab:8cycles} precisely those of rows $9$, $13$ and $17$ hold,
so $7 + v_a(\C_1) = v_a(\C_2) + v_a(\C_3)$ implies $v_a(\C_3) = 10$. 
In view of Lemma~\ref{le:a2n6}, the set $\C_3$ consists of five $H$-orbits of $8$-cycles. 
The existence of each of the $H$-orbits of $8$-cycles of trace $a^2n^6$ is uniquely determined by
a condition of the form 
\begin{equation}
\label{eq:a2n6cycle}
	-1+\delta_1 r + \delta_2 r^2 + \delta_3 r^3 + \delta_4 r^4 + \delta_5 r^5 + t = 0,\quad \delta_i \in \{0,1\},\quad
	i \in \{1,2,3,4,5\}.
\end{equation}
Since $t(r-1)=0$, we either have $t = 0$ or $t = 7$, so the above conditions are easy to check.
It turns out that in case $t=0$ only $-1+r^3+r^4+r^5+t=0$ holds, and so $v_a(\C_3) = 2$, and in case
$t=7$ only conditions $-1+r^2+r^3+t = 0$, $-1+r+r^5+t=0$ and $-1+r+r^2+r^3+r^4+t=0$ hold, so
$v_a(\C_3) = 6$. In both cases $v_a(\C_3) \neq 10$, a contradiction. Thus no $8$-cycle of trace $a^5nan$
exists, as claimed.

We therefore have $v_a(\C) = 2 + v_a(\C_1) + v_a(\C_2) + v_a(\C_3)$
and $v_g(\C) + v_z(\C) = 2 + v_a(\C_1) + 3v_a(\C_2) + 3v_a(\C_3)$. Thus 
$2 + v_a(\C_1) = v_a(\C_2) + v_a(\C_3)$. Since $v_a(\C_1)$ and $v_a(\C_3)$ are both even, so is $v_a(\C_2)$.
Note that $v_a(\C_2) \neq 2$, for otherwise Lemma~\ref{le:possHorbits} implies that
$r^2 = \frac{n}{2} \pm 1$, and so $r^4 = 1$, which, in view of $r^6 = 1$, implies $r^2 = 1$.
Thus $v_a(\C_2) \in \{0,4\}$.

We claim that no $8$-cycle of trace $a^2nan^2an$ exists in $X$ and thus $v_a(\C_1) \leq 4$.
If this is not the case, then, in view of Proposition~\ref{pro:theisomorphisms}, we can assume
that either $1+2r+r^2=0$ or $1+2r-r^2=0$. If $1+2r+r^2=0$, then $0=r+2r^2+r^3 = -1-r+r^2+r^3$, so
the $H$-orbit corresponding to row $25$ of Table~\ref{tab:8cycles} exists, contradicting the fact
that $v_a(\C_2) \in \{0,4\}$. If however $1+2r-r^2=0$, then $r^4 = (1+2r)^2 = 5 + 12r$, so
$r^6 = 29 + 70r$, which, in view of $r^6 = 1$, implies $28 + 70r = 0$. 
Consequently, $70+140r-70r^2 = 0$, so $14 + 28r=0$ and thus also $42+84r=0$, which gives $14+14r=0$. 
Therefore, $n = 14$. But 
then $r=\pm 3$ or $r=\pm 5$, which is impossible in view of the fact that
no $8$-cycle of trace $a^5nan$ exists. Thus no
$8$-cycle of trace $a^2nan^2an$ exists in $X$, as claimed.

Recall that $v_a(\C_2) \in \{0,4\}$. Suppose first that $v_a(\C_2) = 4$. 
Then, in view of Lemma~\ref{le:possHorbits}, we have
$r^3 = \pm 1$. Moreover, in view of $v_a(\C_3) \neq 0$, the equation
$2 + v_a(\C_1) = v_a(\C_2) + v_a(\C_3)$ implies that $v_a(\C_1) = 4$ and $v_a(\C_3) = 2$.
Proposition~\ref{pro:theisomorphisms} thus implies that we can assume 
that one of $2+r+r^2=0$ and $2+r-r^2=0$ holds, which together with $r^3 = \pm 1$ implies that
$8$-cycles of trace $a^2nan^2an$ exist, a contradiction. 
Therefore $\C_2$ is empty and thus $v_a(\C_3) = 2 + v_a(\C_1)$. Since $v_a(\C_1) \leq 4$ is even, we have
three possibilities.

Suppose $v_a(\C_1) = 0$. Then $v_a(\C_3) = 2$ and in addition $v_a(\C) = 2 + 2 = 4$, so
$v_a(\C) = v_g(\C) = v_z(\C)$ implies $v_g(\C_3) = v_z(\C_3)=3$. But $\C_3$ consists of a single
$H$-orbit of $8$-cycles in this case, so since Lemma~\ref{le:a2n6} implies that $v_g(\C_3)$ cannot be odd,
we have a contradiction.

Suppose now that $v_a(\C_1) = 4$ and thus $v_a(\C_3) = 6$. Since $r^3 \neq \pm 1$, 
Lemma~\ref{le:C1orb} implies that we can 
assume $n = 28$ and $r = 5$. Therefore $1+r+r^2+r^3+r^4+r^5=14$, and so $2t=14$, which
forces $t$ to be either $7$ or $21$. However, if $t = 7$, then out of conditions of form $(\ref{eq:a2n6cycle})$
precisely $-1+r^3+r^4+t=0$ and $-1+r+r^5+t=0$ hold and if $t = 21$, only the condition
$-1+r^2+r^3+r^4+r^5+t=0$ holds. Thus $v_a(\C_3) \neq 6$, a contradiction.

We are left with the possibility that $v_a(\C_1) = 2$. 
Then $v_a(\C_3) = 4$, and so $\C_3$ consists of two $H$-orbits
of $8$-cycles. In view of Proposition~\ref{pro:theisomorphisms} we can assume that either
$2+r+r^2=0$ or $2+r-r^2=0$. Suppose $2+r+r^2=0$. Then $v_g(\C) = 1 + 2 + v_g(\C_3)$, so since
$v_g(\C) = v_a(\C) = 2 + 2 + 4 = 8$, we have $v_g(\C_3) = 5$, which is impossible in view of
Lemma~\ref{le:a2n6}.
Therefore, we can now finally assume that $2+r-r^2=0$.
This implies $v_g(\C) = 1 + 1 + v_g(\C_3)$
and $v_z(\C) = 1 +1+v_z(\C_3)$, and so $v_g(\C_3) = v_z(\C_3) = 6$. Let us denote the two
$H$-orbits of $8$-cycles of $X$ that constitute $\C_3$ by $\HH_1$ and $\HH_2$. Since
$v_g(\HH_i)$ is even for $i = 1,2$, we can assume that
either $v_g(\HH_1) = 6$ and $v_g(\HH_2) = 0$ or $v_g(\HH_1) = 4$ and $v_g(\HH_2)=2$.
In the first case the $H$-orbit $\HH_1$ corresponds to condition
$-1+t=0$. But this is impossible in view of $(\ref{eq:Xparamcond})$ and $r \neq 1$.
Thus $v_g(\HH_1) = 4$ and $v_g(\HH_2)=2$.
Note that $2+r-r^2=0$ implies $r^4 = 4+ 4r + r^2 = 6+5r$, so
$r^6 = 22 + 21r$, which forces $21(r+1) = 0$. Suppose now that $\HH_1$ is determined by 
the condition $-1+\delta_1 r + \delta_2 r^2 + \delta_3 r^3 + \delta_4 r^4+ \delta_5 r^5 + t=0$ and
that $\HH_2$ is determined by the condition
$-1+\delta_1' r + \delta_2' r^2 + \delta_3' r^3 + \delta_4' r^4+ \delta_5' r^5 + t=0$
for some $\delta_i, \delta_i ' \in \{0,1\}$, where $i \in \{1,2,3,4,5\}$.
Thus 
\begin{equation}
\label{eq:deltasdiff}
	(\delta_1 - \delta_1')r + (\delta_2 - \delta_2')r^2 + (\delta_3 - \delta_3')r^3 + (\delta_4 - \delta_4')r^4 + (\delta_5 - \delta_5')r^5 = 0. 
\end{equation}
Since $n$ is even and $r \in \ZZ_n^*$, an even number of the numbers 
$\delta_i - \delta_i'$ are nonzero.
Suppose first that four such numbers are nonzero. Then adding the above two equations we get
$$ -2 + r + r^2+ r^3 + r^4 + r^5 + \delta r^i + 2t = 0$$
for some $i \in \{1,2,3,4,5\}$ and some $\delta \in \{-1,1\}$.
Thus $-3 + \delta r^i = 0$. Clearly $i = 1$ and $i = 5$ are impossible in view of the fact that
no $8$-cycle of trace $a^5nan$ exists in $X$. Moreover, $i = 3$ leads to $9 = 1$, which is impossible.
Thus either $i = 2$ or $i = 4$. It follows that $\delta = \delta^3 r^6 = 27$, and so $n = 28$ (for otherwise $n = 14$,
which contradicts the fact that
no $8$-cycle of trace $a^5nan$ exists) or $n = 26$. Since $21(r+1) = 0$, the latter case
is impossible. Therefore, $r^6 = 1$ and $2+r-r^2=0$ imply
$r = 23$ (we now have $n=28$). 
But then $1-r+2r^2=0$, which
contradicts the fact that no $8$-cycle of trace $a^2nan^2an$ exists in $X$.
Therefore, two of the 
numbers $(\delta_i - \delta_i')$ are nonzero. But this implies that multiplying the
equation $(\ref{eq:deltasdiff})$ 
by an appropriate power of $r$ we get $1 + \delta r^i = 0$, for some $i \in \{1,2,3,4,5\}$
and $\delta \in \{-1, 1\}$. Clearly $i$ cannot be $1$, $2$, $4$ or $5$, for otherwise $r^2  = \pm 1$.
But since $\C_2 = \emptyset$, $i = 3$ is also impossible. This shows that
no two conditions for $H$-orbits of $8$-cycles of trace $a^2n^6$ can hold simultaneously in this case,
which completes the proof.
\end{proof}

\begin{lemma}
\label{le:m=6}
Let $X = \X_e(m,n;r,t)$, where $m = 6$, $n \geq 4$ is an even integer and $r\in \ZZ_n^*$, 
$t \in \ZZ_n$ are such that  (\ref{eq:Xparamcond}) holds, and furthermore  $r^2 \neq \pm 1$. 
Then $X$ is half-arc-transitive unless 
$n = 14n_1$, where $n_1$ is coprime to $7$, precisely
one of $\{r, -r, r^{-1}, -r^{-1}\}$ solves the equation $2-x-x^2=0$ and if we let $r'$ be this unique
solution and $t' = t$ in case $r' \in \{r,r^{-1}\}$ and $t' = t + r + r^3 + \cdots + r^{m-1}$ in case
$r' \in \{-r, -r^{-1}\}$, then $r' \equiv 5 \pmod 7$,
$7(r'-1) = 0$, $t' \equiv 0 \pmod 7$ and $2+r'+t'=0$,
in which case it is arc-transitive.
\end{lemma}

\begin{proof}
Let $\C$ denote the set of all $8$-cycles of $X$. Since the group $H$ (recall that $H$ is as in 
Section~\ref{sec:notation})
acts half-arc-transitively on $X$, the graph $X$ can only be half-arc-transitive or arc-transitive.
Suppose it is arc-transitive. Then, in view of Lemma~\ref{le:m=6equal}, 
$v_a(\C)$ differs from at least one of $v_g(\C)$ and $v_z(\C)$. Moreover, Proposition~\ref{pro:halftrcond}
implies that $v_a(\C)$ differs from at most one of $v_g(\C)$ and $v_z(\C)$. In view of 
Proposition~\ref{pro:theisomorphisms} we can therefore assume that $v_a(\C) = v_g(\C)$ and 
$v_a(\C) \neq v_z(\C)$. Note that this implies that for any $\Aut X$-orbit $\C'$ of $8$-cycles of $X$
we have $v_a(\C') = v_g(\C')$ and in addition $\Aut X$ fixes the set $Zig X$.

Let $\C'$ denote the $\Aut X$-orbit of the generic $8$-cycles of $X$. Thus every $8$-cycle of $\C'$
has precisely two zigzags which are antipodal. Table~\ref{tab:8cycles} reveals that,
apart from $8$-cycles of trace $a^2n^6$, only
the $H$-orbits of rows $1$, $10$, $17$, $18$ and $25$ can be in $\C'$ . We claim that some $8$-cycle of trace
$a^2n^6$ exists in $\C'$. Suppose this is not true. Then, in view of Lemma~\ref{le:quartic},
at least one of the $H$-orbits of rows $10$ and $17$ of Table~\ref{tab:8cycles}
is in $\C'$. In view of $v_a(\C') = v_g(\C')$, the $H$-orbit of row $18$ exists and is
contained in $\C'$, and so $r^3 = 1$. Note that this implies that the $H$-orbit of row $25$
cannot exist or else $r = 1$. Therefore, $v_a(\C') = v_g(\C')$ implies that
precisely one of the $H$-orbits of rows $10$ and $17$ of Table~\ref{tab:8cycles} is in $\C'$.
However, $r^3 = 1$ implies that $2-r-r^2 = 0$ if and only if $1+r-2r^2=0$, and so the $H$-orbits 
of rows $10$ and $17$ both exist and precisely one of them is contained in $\C'$. Denote the $\Aut X$-orbit of
the other of the two by $\C''$. Since also $v_a(\C'') = v_g(\C'')$ and since $\C''$ can contain no  
other $H$-orbit corresponding to some row
of Table~\ref{tab:8cycles}, there is some $H$-orbit of $8$-cycles of trace $a^2n^6$ in $\C''$.
But then Lemma~\ref{le:a2n6} implies that $v_a(\C'')$ is even, whereas 
$v_g(\C'')$ is odd, a contradiction.
Therefore $\C'$ contains an $8$-cycle of trace $a^2n^6$, as claimed.

The code of the $8$-cycles of trace $a^2n^6$ that are in $\C'$ is $a^2gzg^3z$. There are precisely two possible
$H$-orbits of $8$-cycles with such codes. The first one, denoted  
by $\HH_1$, is given by the condition $-1+r^2+r^3+r^4+r^5+t=0$ with
a representative $u_0^0u_1^0u_0^{-1}u_1^{-1}u_2^{-1}u_3^{-1+r^2}u_4^{-1+r^2+r^3}u_5^{-1+r^2+r^3+r^4}$
and the second one, denoted by
$\HH_2$, is given by the condition $-1+r+r^2+r^3+r^4+t=0$ with a representative
$u_0^0u_1^0u_0^{-1}u_1^{-1}u_2^{-1+r}u_3^{-1+r+r^2}u_4^{-1+r+r^2+r^3}u_5^{-1+r+r^2+r^3+r^4}$.
In view of $(\ref{eq:Xparamcond})$, the $H$-orbit $\HH_1$ exists if and only if
$2 + r + t = 0$ and $\HH_2$ exists if and only if $2 + r^5 + t = 0$. This implies that
precisely one of $\HH_1$ and $\HH_2$ exists, for otherwise $r = r^5$, which is impossible.
In view of Proposition~\ref{pro:theisomorphisms} some
isomorphism between $\X_e(m,n;r,t)$ and $\X_e(m,n;r^{-1},t)$ preserves the set of anchors, the
set of glides and the set of zigzags, and so we can assume that $\HH_1$ exists. Therefore,
we now have 
\begin{equation}
\label{eq:a2n6orbit}
	-1+r^2+r^3+r^4+r^5+t=0\quad \mathrm{and}\quad 2+r+t=0.
\end{equation}
Since $t(r-1) = 0$, we have $(2+r)(r-1) = 0$, so $2-r-r^2=0$ and thus the
$H$-orbit of row~10 of Table~\ref{tab:8cycles} exists. Suppose that the $H$-orbit
of row~17 of Table~\ref{tab:8cycles} also exists. Then $1+r-2r^2=0$, and so $1-2r+r^2=0$ and
$3-3r=0$, which implies $r^3 = (r-1+1)^3 = 1$. Thus the $H$-orbit of row~18 of Table~\ref{tab:8cycles}
also exists and the $H$-orbit of row~25 does not. But since $v_a(\C') = v_g(\C')$
and $\C'$ contains $\HH_1$, the $H$-orbit of row~18 cannot be in $\C'$. Therefore, if we let $\C''$ denote the
$\Aut X$-orbit of $8$-cycles containing the $H$-orbit of row~18, we have
$v_a(\C'')  < v_g(\C'')$ (which is impossible), 
since the only other $H$-orbit that can be contained in $\C''$ is
one of the $H$-orbits of rows $10$ and $17$ of Table~\ref{tab:8cycles}. Thus $1+r-2r^2 \neq 0$, and
so $\C'$ consists of three $H$-orbits of $8$-cycles, namely the $H$-orbit of generic $8$-cycles, the
one corresponding to row~$10$ of Table~\ref{tab:8cycles} and $\HH_1$. It is now also clear that
none of the $H$-orbits of rows~$18$ and $25$ of Table~\ref{tab:8cycles} can exist.
Since $2-r-r^2=0$ the equation $-1+r^2+r^3+r^4+r^5+t=0$ implies $-5+8r+t=0$, which together with
$(\ref{eq:a2n6orbit})$ implies $7(r-1) = 0$. Thus $7$ divides $n$ in view of $r \neq 1$. Moreover,
$t(r-1) = 0$ implies $t \equiv 0 \pmod 7$, and so $2+r+t=0$ implies $r \equiv 5 \pmod 7$. This also
implies that $7^2$ does not divide $n$ for otherwise $7(r-1) = 0$ implies $r-1 \equiv 0 \pmod 7$.
To summarize
\begin{equation}
\label{eq:7r}
	7(r-1) = 0,\quad 7 \mid n,\quad 7^2 \nmid n, \quad t \equiv 0 \pmod 7,\quad r \equiv 5 \pmod 7 \quad\mathrm{and}\quad 2+r+t=0.
\end{equation}
Recall that we also have $2-r-r^2=0$. Suppose that one of $-r, r^{-1}, -r^{-1}$ would also solve the
equation $2-x-x^2=0$. Clearly it cannot be $-r$. If it was $r^{-1}$, then 
$2r^2-r-1=0$, which is impossible since the $H$-orbit of row~17 of Table~\ref{tab:8cycles} does not
exist. Finally, if it was $-r^{-1}$, then $2r^2 + r - 1=0$, and so $3-r=0$, which is impossible in view of 
$r \equiv 5 \pmod 7$. Thus precisely one of $r,-r,r^{-1},-r^{-1}$ solves the equation
$2-x-x^2=0$. This completes the first part of the proof. 

Let now $n,r$ and $t$ be as in $(\ref{eq:7r})$ and suppose $2-r-r^2=0$. Note that then 
\begin{equation} 
\label{eq:rpow}
	r^2 = 2-r,\  r^3 = -2+3r, \  r^4 = 6 - 5r = -1 + 2r, \  r^5 = 4 - 3r = -3+4r\ \mathrm{and}\ r^5+t=-5+3r.
\end{equation}
We show that $X$ is arc-transitive in this case.
To this end we introduce a certain mapping $\varphi$ of vertices of $X$ and we then
show that it is an automorphism of $X$. The nature of its action will reveal that $X$ is indeed arc-transitive.
The mapping $\varphi$ is defined via Table~\ref{tab:phi}. We define the image of $u_i^j$ under $\varphi$
depending on $i$ and depending on which number $j$ is congruent to modulo $7$.

\begin{table}[htb]
\begin{footnotesize}
\begin{center}
\begin{tabular}{@{}c|c|c|c|c|c|c|c@{}}
	$i \backslash b$ & 0 &  1 & 2 & 3 & 4 & 5 & 6 \\
	\hline  & & & & & & & \\
   $0$ & $u_0^j$ & $u_2^{j+r}$ & $u_2^{j-1+2r}$ & $u_4^{j-2+4r}$ & $u_4^{j+3-r}$ & $u_2^{j+2-r}$ & $u_2^{j+1}$ \\ 
   & & & & & & & \\
   $1$ & $u_1^j$ & $u_1^{j}$ & $u_1^{j-1+r}$ & $u_3^{j+r}$ & $u_5^{j+3-r}$ & $u_3^{j+2-r}$ & $u_1^{j+1-r}$ \\ 
   & & & & & & & \\
   $2$ & $u_0^{j-1}$ & $u_2^j$ & $u_0^{j-2+r}$ & $u_4^{j+r}$ & $u_0^{j+1-2r}$ & $u_2^j$ & $u_0^{j-r}$ \\ 
   & & & & & & & \\
   $3$ & $u_5^{j-3+4r}$ & $u_1^{j-r}$ & $u_1^{j-2+r}$ & $u_5^{j+r}$ & $u_5^{j-1+2r}$ & $u_3^j$ & $u_5^{j-2+3r}$ \\ 
   & & & & & & & \\
   $4$ & $u_0^{j-5+3r}$ & $u_2^{j-r}$ & $u_0^{j-3+r}$ & $u_4^{j+1-r}$ & $u_4^j$ & $u_4^j$ & $u_4^{j-1+r}$ \\ 
   & & & & & & & \\
   $5$ & $u_5^j$ & $u_3^{j-r}$ & $u_1^{j-3+r}$ & $u_3^{j+3-4r}$ & $u_5^j$ & $u_3^{j+2-3r}$ & $u_3^{j+1-2r}$ \\    
   
\end{tabular}
\caption{The entry in $i$-th row and $b$-th column represents the image $u_i^j\varphi$ in case $j \equiv b \pmod 7$.}
\label{tab:phi}
\end{center}
\end{footnotesize}
\end{table}

The mapping $\varphi$ can be obtained by further investigating the structure of $X$ but in order not to
increase the length of the paper even further, we decide to omit this. For our purposes it is enough
to just prove that $\varphi$ is indeed an automorphism of $X$. 
We first show that $\varphi$ is a permutation
of the vertices of $X$. We then show it preserves the adjacency of vertices.

Let us first show that $\varphi$ is bijective. We accomplish this by showing that $\varphi^2$ is the identity. 
For each $i \in \ZZ_6$ and each $j \in \ZZ_n$ we simply
calculate $u_i^j\varphi^2$ using Table~\ref{tab:phi}. For instance 
if $j \equiv 0 \pmod 7$, then $u_0^j\varphi^2 = u_0^j\varphi = u_0^j$.
If $j \equiv 1 \pmod 7$, then $u_0^j\varphi^2 = u_2^{j+r}\varphi = u_0^{j+r-r} = u_0^j$, since
$j+r \equiv 6 \pmod 7$ (recall that $r \equiv 5 \pmod 7$).
If $j \equiv 2 \pmod 7$, then $u_0^j\varphi^2 = u_2^{j-1+2r}\varphi = u_0^{j-1+2r+1-2r} = u_0^j$, since  
$j-1+2r \equiv 4 \pmod 7$. 
Continuing in this fashion one finally finds that $\varphi^2 = \mathrm{Id}$, as desired. We leave the
details to the reader. 

We now show that $\varphi$ preserves the adjacencies of $X$.
We first check that each edge of the form $u_i^j u_{i+1}^j$, where $i \in \ZZ_6 \setminus \{5\}$, $j \in \ZZ_n$,
is mapped to an edge of $X$. Note that we only need to check, that in each column of Table~\ref{tab:phi}
the vertices of two consecutive rows are adjacent. Most of this is obvious in view of $(\ref{eq:rpow})$.
For column labeled $0$, in view of $(\ref{eq:7r})$ and $(\ref{eq:rpow})$, we have 
$j -3 +4r + r^5 + t = j - 1 - 7 + 7r = j-1$, so
$u_0^{j-1}$ is adjacent to $u_5^{j-3+4r}$. Since $t = -2 -r$, the vertex $u_5^{j-3+4r}$ is also adjacent
to $u_0^{j-5+3r}$. Finally, $j + r^5+t=j-5+3r$, which takes care of column $0$. 
We leave the other columns to the reader.

We now check that each edge of the form $u_i^j u_{i+1}^{j+r^i}$, 
where $i \in \ZZ_6 \setminus \{5\}$, $j \in \ZZ_n$, is mapped to an edge of $X$.
Note that $r \equiv 5 \pmod 7$ implies $r^2 \equiv 4 \pmod 7$, $r^3 \equiv 6 \pmod 7$, $r^4 \equiv 2 \pmod 7$ and 
$r^5 \equiv 3 \pmod 7$. In Table~\ref{tab:phi2} 
we list the images of the vertices $u_{i+1}^{j+r^i}$ under $\varphi$. 
Using this table and Table~\ref{tab:phi} our claim is easily
seen to hold as one only needs to check that the vertices at the same entries of the two tables are adjacent.

\begin{table}[htb]
\begin{footnotesize}
\begin{center}
\begin{tabular}{@{}c|l|l|l|l@{}}
	$i \backslash b$ & 0 & 1 & 2 & 3 \\
	\hline  & & & & \\
   $0$ & $u_1^{j+1}$ & $u_1^{j+1-1+r} = u_1^{j+r}$ & $u_3^{j+1+r} = u_3^{j-1+2r+r^2}$ & $u_5^{j+1+3-r} = u_5^{j-2+4r+r^4}$  \\ 
   & & & & \\
   $1$ & $u_2^{j+r}$ & $u_0^{j+r-r} = u_0^{j}$ & $u_0^{j+r-1} = u_0^{j-1+r}$ & $u_2^{j+r}$ \\ 
   & & & & \\
   $2$ & $u_5^{j+r^2-1+2r} = u_5^{j-1-t}$ & $u_3^{j+r^2}$ & $u_5^{j+r^2-2+3r} = u_5^{j-2+r-t}$ & $u_5^{j+r^2-3+4r} = u_5^{j+r+r^4}$ \\ 
   & & & &\\
   $3$ & $u_4^{j+r^3-1+r} = u_4^{j-3+4r}$ & $u_0^{j+r^3-5+3r} = u_0^{j-r}$ & $u_2^{j+r^3-r} = u_2^{j-2+2r}$ & $u_0^{j+r^3-3+r} = u_0^{j+r+r^5+t}$\\ 
   & & & & \\
   $4$ & $u_1^{j+r^4-3+r} = u_1^{j-4+3r}$ & $u_3^{j+r^4+3-4r} = u_3^{j-r+r^2}$ & $u_5^{j+r^4} = u_5^{j-3+r-t}$ & $u_3^{j+r^4+2-3r} = u_3^{j+1-r}$ \\    
\end{tabular}
\vskip 0.5cm

\begin{tabular}{@{}c|l|l|l@{}}
	$i \backslash b$ & 4 & 5 & 6 \\
	\hline  & & & \\
   $0$ & $u_3^{j+1+2-r} = u_3^{j+3-r}$ & $u_1^{j+1+1-r} = u_1^{j+2-r}$ & $u_1^{j+1}$ \\ 
   & & & \\
   $1$ & $u_0^{j+r-2+r} = u_0^{j+3-r+r^5+t}$ & $u_4^{j+r+r} = u_4^{j+2-r+r^3}$ & $u_0^{j+r+1-2r} = u_0^{j+1-r}$ \\ 
   & & & \\
   $2$ & $u_1^{j+r^2-r} = u_1^{j+2-2r}$ & $u_1^{j+r^2-2+r} = u_1^{j}$ & $u_5^{j+r^2+r} = u_5^{j-r-t}$ \\ 
   & & & \\
   $3$ & $u_4^{j+r^3+1-r} = u_4^{j-1+2r}$ & $u_4^{j+r^3}$ & $u_4^{j+r^3} = u_4^{j-2+3r}$ \\ 
   & & & \\
   $4$ & $u_3^{j+r^4+1-2r} = u_3^j$ & $u_5^{j+r^4}$ & $u_3^{j+r^4-r} = u_3^{j-1+r}$\\    
\end{tabular}
\caption{The entry in $i$-th row and $b$-th column represents the image of $u_{i+1}^{j+r^i}$ 
under $\varphi$ in 
case $j \equiv b \pmod 7$.}
\label{tab:phi2}
\end{center}
\end{footnotesize}
\end{table}

Finally, we show that all the edges connecting $X_5$ to $X_0$ are mapped to edges of $X$. This is
demonstrated by Table~\ref{tab:phi3}, where, for each $j$ up to congruence modulo $7$,
the images of vertices $u_0^{j+t}$ and $u_0^{j+r^5+t}$ (which are the neighbors of $u_5^{j}$ in $X_0$) are
given. Recall that $t \equiv 0 \pmod 7$, and so $t + r^5 \equiv 3 \pmod 7$. 
This table thus finally confirms that $\varphi$ is an automorphism of $X$.

\begin{table}[htb]
\begin{footnotesize}
\begin{center}
\begin{tabular}{@{}c|l|l@{}}
	$b \backslash $ & $u_0^{j+t}\varphi$	&  $u_0^{j+r^5+t}\varphi$ \\
	\hline  & & \\
   $0$ & $u_0^{j+t}$ & $u_4^{j+r^5+t-2+4r} = u_4^{j}$  \\ 
   & & \\
   $1$ & $u_2^{j+t+r} = u_2^{j-r-r^2}$ & $u_4^{j+r^5+t+3-r} = u_4^{j-r+r^3}$ \\ 
   & & \\
   $2$ & $u_2^{j+t-1+2r} = u_2^{j-3+r}$ & $u_2^{j+r^5+t+2-r} = u_2^{j-3+2r}$ \\ 
   & & \\
   $3$ & $u_4^{j+t-2+4r} = u_4^{j+3-4r}$ & $u_2^{j+t+r^5+1}=u_2^{j+3-4r}$  \\ 
   & & \\
   $4$ & $u_4^{j+t+3-r} = u_4^{j-r^4}$ & $u_0^{j+r^5+t}$  \\ 
   & & \\
   $5$ & $u_2^{j+t+2-r} = u_2^{j+2-3r-r^2}$ & $u_2^{j+r^5+t+r} = u_2^{j+2-3r}$\\    
   & & \\
   $6$ & $u_2^{j+t+1} = u_2^{j+1-2r - r^2}$ & $u_2^{j+r^5+t-1+2r} = u_2^{j+1-2r}$
   
\end{tabular}
\caption{The entries in $b$-th row represent the images of $u_0^{j+t}$ and $u_0^{j+r^5+t}$ under $\varphi$, 
respectively, in case $j \equiv b \pmod 7$.}
\label{tab:phi3}
\end{center}
\end{footnotesize}
\end{table}

However, if we let $\psi$ denote the unique nontrivial automorphism of $H_{u_1^0}$, then
$u_1^0\varphi\psi\sigma = u_1^0\psi\sigma = u_2^0$ and $u_2^0\varphi\psi\sigma = u_0^{-1}\psi\sigma =
u_0^0\sigma=u_1^0$, and so $\varphi\psi\sigma$ interchanges adjacent vertices $u_1^0$ and $u_2^0$
and thus $X$ is arc-transitive, as claimed.
\end{proof}


\section{The case $m = 4$}
\label{sec:m=4}

\indent 
Throughout this section we let $m=4$ and we let $n,r,t,X,H,\rho,\sigma$ and $\tau$ be as in Section~\ref{sec:notation}.
Recall that this implies that the conditions (\ref{eq:Xparamcond}) are satisfied. In view of Proposition~\ref{pro:r^2}
we can assume $r^2 \neq \pm 1$. Recall also that we can assume 
(\ref{eq:2r}) and (\ref{eq:n}). In view of the results of
Section~\ref{sec:8cycles}, the only possible $H$-orbits of $8$-cycles of $X$ are those of Table~\ref{tab:8cycles}
and the ones with traces $a^4n^4$, $a^2na^2n^3$, $a^2n^2a^2n^2$, $anan^5$ and $n^8$. The next two 
lemmas give some information about these $H$-orbits.

\begin{lemma}
\label{le:special8cycles}
Let $X = \X_e(m,n;r,t)$, where $m = 4$, $n \geq 4$ is an even integer and $r\in \ZZ_n^*$, 
$t \in \ZZ_n$ are such that  (\ref{eq:Xparamcond}) holds, and furthermore  $r^2 \neq \pm 1$. 
Then the following hold.
\begin{itemize}
	\item[(i)] An $8$-cycle of trace $a^4n^4$ exists in $X$ if and only if, up to isomorphisms
	of Proposition~\ref{pro:theisomorphisms}, $X$ is one of 
	$\X_e(4,20;7,10)$ and $\X_e(4,30;13,25)$. 
	In this case precisely one $H$-orbit of $8$-cycles of trace $a^4n^4$ exists in $X$.
	\item[(ii)] An $8$-cycle of trace $a^2na^2n^3$ exists in $X$ if and only if, 
	up to isomorphisms of Proposition~\ref{pro:theisomorphisms}, $X = \X_e(4,30;13,25)$. In this case
	precisely one $H$-orbit of \hbox{$8$-cycles} of trace $a^2na^2n^3$ exists in $X$.
	\item[(iii)] At most one $H$-orbit of $8$-cycles of trace $a^2n^2a^2n^2$ exists in $X$ unless an
	$8$-cycle of trace $a^5nan$ exists in $X$.
	\item[(iv)] Each $H$-orbit of $8$-cycles of trace $a^4n^4$, $a^2na^2n^3$ or $anan^5$ has length $8n$ and
	each $H$-orbit of $8$-cycles of trace $a^2n^2a^2n^2$ has length $4n$.	
\end{itemize}
\end{lemma}

\begin{proof}
Suppose that an $8$-cycle of trace $a^4n^4$ exists.
Note that each $H$-orbit of $8$-cycles of trace $a^4n^4$
is uniquely determined by a condition of the form 
$$ 3+\delta_1 r + \delta_2 r^2 + \delta_3 r^3 + t = 0,\ \mathrm{where}\ \delta_1, \delta_2, \delta_3 \in \{0,1\}.$$ 
Since our consideration is only modulo isomorphisms of Proposition~\ref{pro:theisomorphisms}, we only need to
consider the conditions $3+t=0$, $3+r+t=0$, $3+r^2+t=0$ and $3+r+r^2+t=0$. 
If $3+t=0$, then (\ref{eq:Xparamcond}) implies $3(r-1) = 0$, and so $3r=3$. Since $2t=-6$ and
$1+r+r^2+r^3 +2t = 0$, we have $0 = -15 + 3r+3r^2+3r^3 =-6$, which contradicts $(\ref{eq:n})$. 
If $3+r^2+t=0$, then $3+r^2+3r+r^3+2t=0$, and so $2+2r=0$, which
contradicts (\ref{eq:2r}). 
Suppose now that $3+r+t=0$. Then $(3+r)(r-1)=r^2+2r-3=0$. Also $3+r+3r^2+r^3+2t=0$, so $2+2r^2=0$ and thus,
in view of $r^2\neq \pm 1$, we have $r^2 = \frac{n}{2}-1$. Therefore, $2r^2+4r-6=0$ implies $4r-8=0$. Moreover,
$1+r+r^2+r^3 = \frac{n}{2}(1+r) = 0$, and so $2t=0$. In view of $3+r+t=0$, we have $12+4r=0$ and thus 
$n=20$. Then $6+2r=0$ implies that either $r = 7$ or $r = 17$. In the former case
$3+r+t=0$ implies $t = 10$ and in the latter case $t = 0$. Thus either
$X = \X_e(4,20;7,10)$ or $X = \X_e(4,20;17,0)$ (note that $17 = -7^{-1}$ in $\ZZ_{20}$). 
For both of these graphs $3+r+t=0$, so an $H$-orbit of
$8$-cycles of trace $a^4n^4$ does exist. Suppose finally that $3+r+r^2+t=0$.
Then (\ref{eq:Xparamcond}) implies $2-r^3-t=0$, and so multiplication by $-r$ gives $1-2r+t=0$.
Subtracting this from $3+r+r^2+t=0$ we get $2+3r+r^2=0$. We also have $3r+r^2+r^3+t=0$, so
$3-2r-r^3=0$ and thus also $1-3r+2r^2=0$. Together with $2+3r+r^2=0$, this implies $3+3r^2=0$,
and so $6+9r+3r^2=0$ gives $3+9r=0$. Moreover, $9 - 6r - 3r^3=0$ implies $9-3r=0$. Thus $27-9r=0$, and
so $n=30$. Then $3r=9$ implies (recall that $r \in \ZZ_{30}^*$) that $r$ is either $13$ or $23$.
In the former case $t = 25$ and in the latter case $t = 15$. Therefore $X$ is one of
$\X_e(4,30;13,25)$ and $\X_e(4,30;23,15)$ (note that $23 = -13^{-1}$ in $\ZZ_{30}$). 
For these two graphs $3+r+r^2+t=0$, so
an $H$-orbit of $8$-cycles of trace $a^4n^4$ does exists in $X$. Clearly no nontrivial
element of $H$ can fix an $8$-cycle of trace $a^4n^4$, and so all $H$-orbits of $8$-cycles of trace $a^4n^4$
have length $8n$, as claimed. The arguments of this paragraph also imply that $X$ cannot 
have two $H$-orbits of $8$-cycles of trace $a^4n^4$.

Suppose now an $8$-cycle of trace $a^2na^2n^3$ exists.
Note that each $H$-orbit of $8$-cycles of trace $a^2na^2n^3$
is uniquely determined by a condition of the form 
$$ 2+\delta_1 r + \delta_2 r^2 + \delta_3 r^3 + t = 0,\ \mathrm{where}\ 
\delta_1 \in \{-1,2\}\ \mathrm{and}\ \delta_2, \delta_3 \in \{0,1\}.$$ 
Using Proposition~\ref{pro:theisomorphisms}, one can check that we only need to consider
the conditions $2-r+t=0$, $2-r+r^2+t=0$ and $2-r+r^3+t=0$. 
Suppose first that
$2-r+t=0$. It follows that $(2-r)(r-1)=0$, and so $2-3r+r^2=0$. Then $r^4 = (-2+3r)^2 = 4-12r+9r^2 = -14+15r$, and
so $r^4 = 1$ implies $15(r-1)=0$. On the other hand $r^3 = -2r+3r^2=-6+7r$ implies
$1+r+r^2+r^3+2t=-7+11r+2t=0$, and so $2t=-4+2r$ (in view of $2-r+t=0$) 
implies $-11+13r=0$. Thus $4-2r=0$, and so $13+r=0$.
Hence $26+2r=0$, and so $n = 30$. Therefore, $r = 17$ and consequently $t= 15$. Thus $X = \X_e(4,30;17,15)$
which is isomorphic with $\X_e(4,30;23,15)$ via an isomorphism of Proposition~\ref{pro:theisomorphisms}.
Note that in this graph an $H$-orbit of $8$-cycles of trace $a^2na^2n^3$ does exists.
Suppose now that $2-r+r^2+t=0$. Then $2r-r^2+r^3+t=0$, and so $2+r+r^3+2t=0$, forcing $1-r^2=0$, a contradiction.
Suppose finally
that $2-r+r^3+t=0$. Then $2r^2-r^3+r+t=0$ and thus $2+2r^2+2t=0$, which implies $1-r+r^2-r^3=0$.
On the other hand we have $2r-r^2+1+t=0$, so subtracting this equation from $2-r+r^3+t=0$ 
yields $1-3r+r^2+r^3=0$. This in turn implies $2r-2r^3=0$, and so
$2-2r^2=0$, forcing $r^2 = \frac{n}{2} + 1$. But then $0 = 1-r+r^2-r^3 = 1 + \frac{n}{2}+1-r(1 + \frac{n}{2} + 1)=
2-2r$, which contradicts $(\ref{eq:2r})$. Clearly no element of $H$ can fix an $8$-cycle of trace
$a^2na^2n^3$, and so each such $H$-orbit is indeed of length $8n$, as claimed.
The arguments of this paragraph also imply that $X$ cannot 
have two $H$-orbits of $8$-cycles of trace $a^2na^2n^3$.

Observe that an $H$-orbit of $8$-cycles of trace $a^2n^2a^2n^2$
is uniquely determined by a condition of the form 
$$2+\delta_1 r + \delta_2 r^2 + \delta_3 r^3 + t = 0,\ \mathrm{where}\ 
\delta_2 \in \{-1,2\}\ \mathrm{and}\ \delta_1, \delta_3 \in \{0,1\}.$$ 
Since $r^2 \neq \pm 1$, it is clear that
no two conditions of the form $2+\delta_1r - r^2 + \delta_3 r^3+t=0$ nor of the form
$2+\delta_1 r + 2r^2 + \delta_3r^3 + t=0$ can hold simultaneously. Suppose 
that $2+\delta_1 r - r^2 + \delta_3 r^3 + t= 0 = 2+\delta_1'r+2r^2+\delta_3'r^3+t$ for some
$\delta_1,\delta_1',\delta_3,\delta_3' \in \{0,1\}$. Then 
$(\delta_1 - \delta_1')r -3r^2+ (\delta_3 - \delta_3')r^3 = 0$. Since $n$ is even and $r \in \ZZ_n^*$,
precisely one of $(\delta_1 - \delta_1')$
and $(\delta_3 - \delta_3')$ is zero. (The other is thus $\pm 1$.) Thus either $-3r^2 \pm r^3 = 0$ and 
consequently $-3 \pm r = 0$, or $\pm r - 3r^2=0$ and 
consequently $\pm 1 - 3r = 0$, and so an $8$-cycle of trace $a^5nan$ exists.
We now show that each $H$-orbit of $8$-cycles of trace $a^2n^2a^2n^2$ is of length $4n$. Let $C$ be any
$8$-cycle of trace $a^2n^2a^2n^2$, such that
$u_1^1u_0^0u_1^0$ is one of its anchors. This is the only anchor of $C$, whose internal vertex is in $X_0$. 
Moreover, $C$ also
contains precisely one nonanchor, whose internal vertex is in $X_0$, and it does not contain any negative anchor, whose
internal vertex is in $X_1$. It is therefore clear, that
no element of $\langle \rho, \tau \rangle$
or $ \sigma\langle \rho, \tau\rangle$ can fix $C$ setwise and thus the $H$-orbit of $C$ has length at least $4n$.
We claim that an element of $ \sigma^2\langle \rho, \tau\rangle$ does fix $C$, however, proving that
the length is indeed $4n$. Suppose on the contrary that this is not the case. We can assume that $C$ also
contains $u_0^{-1}$ (otherwise $C\tau$ does). It is easy to see that, 
up to isomorphisms of Proposition~\ref{pro:theisomorphisms},
the corresponding condition is then one of $2-r^2+t=0$ and $2+2r^2+r^3+t=0$. 
However, these conditions are both contradictory. Namely, if $2-r^2+t=0$, then
$(\ref{eq:Xparamcond})$ implies $1-r-2r^2-r^3-t=0$, so also $-2-r+r^2-r^3-t=0$ and thus $-r-r^3=0$, which
contradicts $r^2 \neq -1$. As for $2+r+2r^2+t=0$, note that in this case
$2r^2+r^3+2+t=0$, so subtracting the two equations we get $r^3-r=0$, which is also impossible. 

Finally, it is clear that no nontrivial element of $H$ can fix an $8$-cycle of trace $anan^5$. This completes
the proof.
\end{proof}

\begin{lemma}
\label{le:4coiled}
Let $X = \X_e(m,n;r,t)$, where $m = 4$, $n \geq 4$ is an even integer and $r\in \ZZ_n^*$, 
$t \in \ZZ_n$ are such that  (\ref{eq:Xparamcond}) holds, and furthermore  $r^2 \neq \pm 1$. 
Then the possible $H$-orbits of coiled $8$-cycles of $X$ are 
those represented in 
Table~\ref{tab:4coiled}, where a representative, the code, a necessary and sufficient condition
for its existence and the length of the $H$-orbit are given. 
Moreover, the $H$-orbit corresponding to row~$4$ cannot exist
simultaneously with any of the $H$-orbits corresponding to rows $1$ or $5$.
\end{lemma}

\begin{table}[!hbt]
\begin{footnotesize}
\begin{center}
\begin{tabular}{@{}c|l|c|c|c@{}}
	Row	&  A representative 	& Code	& Condition	& Orbit length \\
	\hline  & & & & \\
   $1$ & $u_0^0u_1^0u_2^0u_3^0u_0^tu_1^tu_2^tu_3^t$ & $g^8$ & $2t = 0\ \mathrm{and}\ t \neq 0$ & $n$ \\    
   $2$ & $u_0^0u_1^0u_2^0u_3^0u_0^tu_1^{1+t}u_2^{1+r+t}u_3^{1+r+r^2+t}$ & $g^3zg^3z$ & $t \neq 0$ & $4n$ \\
   $3$ & $u_0^0u_1^0u_2^0u_3^{r^2}u_0^{r^2+t}u_1^{1+r^2+t}u_2^{1+r+r^2+t}u_3^{1+r+r^2+t}$ & $gz^3gz^3$ & $1+r^2+t \neq 0$ & $4n$ \\
   $4$ & $u_0^0u_1^0u_2^0u_3^{r^2}u_0^{r^2+r^3+t}u_1^{r^2+r^3+t}u_2^{r^2+r^3+t}u_3^{2r^2+r^3+t}$ & $gzgzgzgz$ & $2+2r+2t=0$ & $2n$ \\
   $5$ & $u_0^0u_1^0u_2^ru_3^ru_0^{r+r^3+t}u_1^{r+r^3+t}u_2^{2r+r^3+t}u_3^{2r+r^3+t}$ & $z^8$ & $2+2r^2+2t=0 \ \mathrm{and}\ 1+r^2+t\neq 0$ & $n$ \\   
\end{tabular}
\caption{Possible $H$-orbits of coiled $8$-cycles of $X$ when $m=4$.}
\label{tab:4coiled}
\end{center}
\end{footnotesize}
\end{table}

\begin{proof}
Note that the existence of a coiled $8$-cycle implies that a condition of the form 
$$ \delta_0 + \delta_1 r + \delta_2 r^2 + \delta_3 r^3 + 2t = 0,\ 
	\mathrm{where}\ \delta_i \in \{0,1,2\}\ \mathrm{for}\ \mathrm{all}\ i \in \{0,1,2,3\}$$ 
holds. Let $s = \delta_0 + \delta_1 + \delta_2 + \delta_3$. 
We can assume $s \leq 4$, since we only need to determine which $H$-orbits of the coiled $8$-cycles
are possible. Moreover, since $n$ is even and $r\in\ZZ_n^*$, it is clear that
$s$ is even. We distinguish three possible cases depending on $s$.

Suppose $s = 0$. Then the condition is $2t=0$ and the corresponding $H$-orbit is clearly
the one represented in row~1 of Table~\ref{tab:4coiled}. Note that the condition $t \neq 0$ is necessary
to ensure that we indeed have an $8$-cycle. The orbit length is also clear. 

Suppose $s = 2$. We can thus assume that the condition is one of
$2+2t=0$, $1+r+2t=0$ and $1+r^2+2t=0$.
Note however that none of them is possible. Namely, multiplying by $r-1$, 
the first implies $2(r-1)=0$, the second implies 
$r^2-1=0$ and the third implies (since $1+r+r^2+r^3+2t=0$) $r+r^3=0$, and so $1+r^2=0$.

Suppose finally $s = 4$. If none of $\delta_i$, where $i \in \{0,1,2,3\}$, is $2$, then
$1+r+r^2+r^3+2t=0$, which always holds. It is easy to see that the only possible $H$-orbits are then
the ones corresponding to rows $2$ and $3$ of Table~\ref{tab:4coiled}. Moreover, for the given representatives
to actually be $8$-cycles the conditions listed are needed. 
We can now assume that $\delta_0 = 2$. Then
precisely one of $\delta_1$, $\delta_2$, $\delta_3$ is $2$ and the other two are zero or else 
$1 - r^i=0$ for some $i \in \{1,2,3\}$, which is clearly impossible. Thus the two possible
conditions are $2+2r+2t=0$ and $2+2r^2+2t=0$. It is clear that each of the two conditions
uniquely determines the corresponding $H$-orbit. 
Moreover, if $2+2r^2+2t=0$ is to give rise to an $H$-orbit of
coiled $8$-cycles, then clearly the condition $1+r^2+t\neq 0$ must hold. The details are given in
rows $4$ and $5$ of Table~\ref{tab:4coiled}. It is easy to check that the lengths of the orbits are indeed
as listed in the table.

The last claim of the lemma is straightforward. 
Namely, if $2+2r+2t=0$, then any of the conditions of rows $1$ and $5$ of Table~\ref{tab:4coiled} contradicts
$(\ref{eq:2r})$.
\end{proof}

\begin{lemma}
\label{le:m=4_eq}
Let $X = \X_e(m,n;r,t)$, where $m = 4$, $n \geq 4$ is an even integer and $r\in \ZZ_n^*$, 
$t \in \ZZ_n$ are such that  (\ref{eq:Xparamcond}) holds, and furthermore  $r^2 \neq \pm 1$. Suppose that 
$\Aut X$ does not act transitively on the set of $2$-paths of $X$. Then
$X$ is half-arc-transitive.
\end{lemma}

\begin{proof}
Suppose on the contrary that $X$ is arc-transitive.
Note that, in view of Proposition~\ref{pro:arc-tr} and Proposition~\ref{pro:theisomorphisms}, 
we can assume that $\Aut X$ has precisely two orbits on the set of $2$-paths, and that one orbit is the set of
anchors and glides and the other is the set of zigzags of $X$. 
Moreover, for any $\Aut X$-orbit $\C'$ of $8$-cycles of $X$ we have $v_a(\C') = v_g(\C')$.

Let $\C$ denote the set of all $8$-cycles of $X$ having precisely two zigzags which are antipodal and 
let $\C'$ denote the $\Aut X$-orbit of $8$-cycles containing the generic $8$-cycles.
Clearly $\C' \subseteq \C$.
Note that $\C$ can consist only of the generic $8$-cycles, of $H$-orbits of $8$-cycles 
corresponding to rows $10$, $17$, $18$ and $25$
of Table~\ref{tab:8cycles} and of $H$-orbits of $8$-cycles of codes $a^2za^2gzg$, $a^2zga^2zg$, 
$azag^2zg^2$, $agazg^3z$ and $zg^3zg^3$. However, the $H$-orbit of row~$18$ of Table~\ref{tab:8cycles} 
cannot exist since $r^3 \neq 1$. Observe that despite the fact that $m = 4$, 
the only condition for the
$H$-orbit of row $25$ to exist is the one stated in Table~\ref{tab:8cycles} since the $8$ vertices
of a representative are indeed distinct (otherwise $r^2+r^3+t=0$, but then multiplication by $(r-1)$ gives $r^2 = 1$,
which is impossible).
Moreover, the proof of Lemma~\ref{le:special8cycles} shows
that no $8$-cycle of code $a^2za^2gzg$ exists as otherwise $2-r+r^3+t=0$, which we saw was impossible.
Furthermore, at most one of the $H$-orbits corresponding to rows $10$ and $17$ of Table~\ref{tab:8cycles}
can exist since otherwise $1-2r+r^2=0$ and $3-3r=0$, and thus
$r^3 = (r-1+1)^3 = 1$, a contradiction. Note also that at most one $H$-orbit of $8$-cycles of $X$
of code $a^2zga^2zg$ can exist. Otherwise, the proof of Lemma~\ref{le:special8cycles} reveals that 
$2-r^2+r^3+t=0$ and $2+r-r^2+t=0$. But then $r^3-r=0$, a contradiction. 
Let $\C_1$ denote the set of $8$-cycles 
corresponding to rows $10$ and $17$ of Table~\ref{tab:8cycles}, let $\C_2$ denote the set of $8$-cycles 
corresponding to row~25 of Table~\ref{tab:8cycles}, let
$\C_3$ denote the set of $8$-cycles 
of code $a^2zga^2zg$, let $\C_4$ denote the set of $8$-cycles of codes $azag^2zg^2$ and $agazg^3z$
and let $\C_c$ denote the coiled $8$-cycles of code $g^3zg^3z$. Note that the above
remarks imply that each of $\C_1$ and $\C_3$ consists of at most one $H$-orbit of $8$-cycles of $X$.
Moreover, $2v_a(\C_2) = v_g(\C_2)$ and $2v_a(\C_4) = v_g(\C_4)$.

Suppose first that $\C_c$ is nonempty. We claim that then $\C_c \subset \C'$. If this is not the case, 
then denote the $\Aut X$-orbit containing $\C_c$ by $\C''$. Since 
$\C'' \subset \C$, the remarks of the previous paragraph show that then $v_a(\C'') < v_g(\C'')$,
a contradiction. Thus $\C_c \subset \C'$. But then $v_a(\C') = v_g(\C')$ and the remarks of the
previous paragraph imply that $\C_1$ and $\C_3$ each contain precisely one $H$-orbit, and that
$\C' = \G \cup \C_1 \cup \C_3 \cup \C_c$, where $\G$ is the set of the generic $8$-cycles of $X$.
Moreover, $\C_2$ and $\C_4$ are both empty. (Otherwise if we let $\C''$ be an $\Aut X$-orbit containing some 
of these $8$-cycles we have $v_a(\C'') < v_g(\C'')$.) Since $r^4=1$, 
the condition $1+r-2r^2=0$ implies $2-r^{-1}-(r^{-1})^2=0$. 
In view of Proposition~\ref{pro:theisomorphisms}
the isomorphism between $\X_e(m,n;r,t)$ and $\X_e(m,n;r^{-1},t)$ preserves the sets of glides, anchors and zigzags and
therefore we can assume that the $H$-orbit corresponding to row~10 of Table~\ref{tab:8cycles} exists, that is 
$2-r-r^2=0$. Note that then $2+r-r^2+t=0$ cannot hold for otherwise $2r+t=0$, 
and so $2r(r-1)=0$, which contradicts $(\ref{eq:2r})$.
Thus $\C_3$ is the $H$-orbit corresponding to condition $2-r^2+r^3+t=0$ and its representative
is $u_0^0u_1^1u_0^1u_1^2u_2^2u_3^2u_2^{2-r^2}u_3^{2-r^2}$. As $X$ was assumed to be arc-transitive,
Proposition~\ref{pro:arc-tr} implies that there exists an automorphism $\varphi$ of $X$ mapping the
anchor $u_1^1u_0^0u_1^0$ to the glide $u_3^{-t}u_0^0u_1^0$. We can in fact assume that
$u_1^0\varphi = u_1^0$, otherwise take $\tau\varphi$. Since $\Aut X$ fixes the set $Zig X$, the zigzag
$u_0^0u_1^0u_2^r$ is mapped to a zigzag, and so $u_2^r\varphi = u_2^r$. There are precisely three
$8$-cycles of $\C'$ containing the $3$-path $u_1^1u_0^0u_1^0u_2^r$, namely one generic $8$-cycle,
one from the set $\C_1$ and one from the set $\C_3$. However, its image under $\varphi$, the $3$-path
$u_3^{-t}u_0^0u_1^0u_2^r$ is contained on precisely one $8$-cycle of $\C'$, namely on a coiled $8$-cycle.
In view of Proposition~\ref{pro:equalfreq}, this contradicts the fact that $\varphi$ is an automorphism of $X$. 

Suppose now that $\C_c$ is empty. Lemma~\ref{le:4coiled} implies that then $t=0$. But then 
$u_0^0u_1^0u_2^0u_3^0$ is a $4$-cycle consisting of four glides. Since $X$ was assumed to be arc-transitive,
Proposition~\ref{pro:arc-tr} implies that there exists an automorphism $\varphi$ of $X$ mapping the
anchor $u_1^1u_0^0u_1^0$ to the glide $u_3^0u_0^0u_1^0$. But then there exists a $4$-cycle of $X$
containing an anchor, and so either $2=0$ or $r \pm 1 = 0$, which are both impossible.
This contradiction shows that $X$ is half-arc-transitive, as claimed.
\end{proof}

\begin{lemma}
\label{le:m=4}
Let $X = \X_e(m,n;r,t)$, where $m = 4$, $n \geq 4$ is an even integer and $r\in \ZZ_n^*$, 
$t \in \ZZ_n$ are such that  (\ref{eq:Xparamcond}) holds, and furthermore  $r^2 \neq \pm 1$. 
Then $X$ is half-arc-transitive.
\end{lemma}

\begin{proof}
In view of Lemma~\ref{le:m=4_eq} we can assume that $\Aut X$ acts transitively on the set of all
$2$-paths of $X$. This implies that the girth of $X$ cannot be $4$ since otherwise $X$ has $4$-cycles
containing anchors which is impossible as the proof of Lemma~\ref{le:m=4_eq} shows.
Let $\C$ be the set of all $8$-cycles of $X$. Then Proposition~\ref{pro:equalfreq}
implies $v_a(\C) = v_g(\C) = v_z(\C)$.
We let $\C_1$ denote the set of $8$-cycles of traces $a^2nan^2an$ and $a^3n^2an^2$, we let
$\C_2$ denote the set of $8$-cycles of trace $an^3an^3$, we let $\C_3$ denote the set of
$8$-cycles of traces $a^4n^4$, $a^2na^2n^3$, $a^2n^2a^2n^2$ and $anan^5$ and we
let $\C_c$ denote the set of coiled $8$-cycles of $X$. 

Clearly $t \neq 0$ and $1+r^2+t \neq 0$, as otherwise $4$-cycles exist in $X$.
Thus Lemma~\ref{le:4coiled} implies that $\C_c$ contains at least the
$H$-orbits corresponding to rows $2$ and $3$ of Table~\ref{tab:4coiled}. Moreover, 
$4 \leq v_g(\C_c), v_z(\C_c) \leq 5$.

We claim that no $8$-cycle of trace $a^5nan$ exists in $X$. If this is not the case, then,
in view of Proposition~\ref{pro:theisomorphisms}, we can assume $r =3$. Thus
$t(r-1) = 0$ implies $2t=0$, and so $0 = 1+r+r^2+r^3 + 2t = 40$, which forces $n$ to
be either $20$ or $40$. If 
$n = 20$, then $t = 10$, since $t \neq 0$. But then $1+r^2+t=0$, which is also impossible.
Thus $n = 40$ and $t = 20$. Therefore $X = \X_e(4,40;3,20)$. Note that
Lemma~\ref{le:special8cycles} implies that there are no $8$-cycles of traces 
$a^4n^4$ or $a^2na^2n^3$ in $X$. The reader can check that 
precisely two $H$-orbits of $8$-cycles of $X$ of trace $a^2n^2a^2n^2$ exist, namely one
corresponding to condition $2-r^2+r^3+t=0$ (and thus of code $a^2zga^2zg$) 
and one to condition $2+2r^2+t=0$ (and thus of code $a^2z^2a^2z^2$), and no
$8$-cycles of trace $anan^5$ exist. Moreover, Table~\ref{tab:4coiled} reveals that 
$\C_c$ consists precisely of the $H$-orbits corresponding to rows $1$, $2$ and $3$.
It is easy to see that $\C_1$ is empty whereas $\C_2$ consists of the $H$-orbit corresponding to
row~22 of Table~\ref{tab:8cycles}. Thus 
$v_a(\C) = 10$ and 
$v_z(\C) =  9$, which contradicts the
fact that $v_a(\C) = v_z(\C)$. Therefore, no $8$-cycle of trace $a^5nan$ exists, as claimed.

Note that this fact and Lemma~\ref{le:C1orb} imply that $\C_1$ consists of at most one $H$-orbit.
Namely, since $r^4 = 1$ but $r^2 \neq 1$, $r$ is of order $4$ in $\ZZ_n^*$.
Suppose an $8$-cycle of trace $a^4n^4$ exists in $X$. In view of Lemma~\ref{le:special8cycles} and the
fact that no $8$-cycle of trace $a^5nan$ exists in $X$, we can assume (Proposition~\ref{pro:theisomorphisms}),
that $X = \X_e(4,30;7,25)$. It is easy to check that then $\C_1$ consists 
of the $H$-orbit corresponding to row~17 of Table~\ref{tab:8cycles}, that
is $1+r-2r^2=0$, and that $\C_2$ consists 
of the $H$-orbit corresponding to row~27 of that table, that is $1-r+r^2-r^3=0$. 
Moreover, $\C_c$ consists of the
$H$-orbits corresponding to rows $2$, $3$ and $5$ of Table~\ref{tab:4coiled}. Furthermore, 
Lemma~\ref{le:special8cycles} implies that precisely one $H$-orbit of $8$-cycles of trace
$a^4n^4$ exists, namely the one corresponding to condition $3+r^2+r^3+t=0$, and thus of code $a^4z^2g^2$, 
and that precisely one
$H$-orbit of $8$-cycles of trace $a^2na^2n^3$ exists, namely the one corresponding to condition
$2+2r+r^2+t=0$ and thus of code $a^2ga^2gz^2$. 
It is also easy to check that no $8$-cycle of trace $a^2n^2a^2n^2$ exists and that precisely
one $H$-orbit of $8$-cycles of trace $anan^5$ exists, namely the one containing
the $8$-cycle $u_0^0u_1^1u_0^1u_3^{1-t}u_0^{1+r^3}u_1^{2+r^3}u_2^{2+r+r^3}u_3^{2+r+r^3}$, which is
of code $agag^2z^2g$.
Thus $v_a(\C) = 15$ and $v_g(\C) = 14$,  
which contradicts $v_a(\C) = v_g(\C) = v_z(\C)$.

Therefore, no $8$-cycle of trace $a^4n^4$ and hence by Lemma~\ref{le:special8cycles} also of
trace $a^2na^2n^3$ exists. Lemma~\ref{le:special8cycles} also implies that at most one $H$-orbit
of $8$-cycles of trace $a^2n^2a^2n^2$ exists. Observe that, in view of the fact that
$\C_1$ contains at most one $H$-orbit of $8$-cycles, we have $2v_a(\C_1) - (v_g(\C_1) + v_z(\C_1)) \leq 4$.
Corollary~\ref{cor:summs} implies 
$2v_a(\C_2) - (v_g(\C_2) + v_z(\C_2)) \leq 0$ and in view of Lemma~\ref{le:4coiled}, we also have
$2v_a(\C_c) - (v_g(\C_c) + v_z(\C_c)) \leq -8$.
Let $\C_{3,1}$ denote the $8$-cycles of $\C_3$
of trace $a^2n^2a^2n^2$ and let $\C_{3,2}$ denote the $8$-cycles of $\C_3$ of trace $anan^5$.
Recall that $\C_{3,1}$ consists of at most one $H$-orbit of $8$-cycles and that, in view of 
Lemma~\ref{le:special8cycles}, such $H$-orbit is of length $4n$.
Therefore $2v_a(\C_{3,1}) - (v_g(\C_{3,1}) + v_z(\C_{3,1})) \leq 2$ and
$2v_a(\C_{3,2}) - (v_g(\C_{3,2}) + v_z(\C_{3,2})) \leq 0$.
Since $v_a(\C) = v_g(\C) = v_z(\C)$, we have $2v_a(\C) - (v_g(\C) + v_z(\C)) = 0$, and so
$$ \begin{array}{ccl}
	0  & =  & 4 - 2 + 2v_a(\C_1) - (v_g(\C_1) + v_z(\C_1)) + 
	2v_a(\C_2) - (v_g(\C_2) + v_z(\C_2)) + \\
	& & 2v_a(\C_{3,1}) - (v_g(\C_{3,1}) + v_z(\C_{3,1})) + 2v_a(\C_{3,2}) - (v_g(\C_{3,2}) + v_z(\C_{3,2})) + \\
	& & 2v_a(\C_c) - (v_g(\C_c) + v_z(\C_c)) \\
	& \leq & 2 + 4 + 0 + 2 + 0 -8 = 0. \end{array}$$
Therefore, $2v_a(\C_1) - (v_g(\C_1) + v_z(\C_1)) = 4$ and thus $8$-cycles of trace $a^2nan^2an$ exist. 
Moreover, $\C_2 = \emptyset$ (otherwise $2v_a(\C_2) - (v_g(\C_2) + v_z(\C_2)) < 0$).
In view of Proposition~\ref{pro:theisomorphisms} we can assume that either $1+2r-r^2=0$ or 
$1+2r+r^2 = 0$. In the former case we have $r^4 = (1+2r)^2 = 5+12r$, and so $r^4 = 1$ implies $4+12r=0$. But then,
in view of $12+24r-12r^2=0$, we have $4+4r=0$, and so $12+12r=0$, which forces $8 \equiv 0 \pmod n$, 
contradicting $(\ref{eq:n})$.
Thus $1+2r+r^2=0$ and hence $r^4 = -3-4r$, so $r^4 = 1$ implies $4(1+r) = 0$. Since $2(1+r) \neq 0$, we
have $2(1+r) = \frac{n}{2}$. Then $r^2 = (r+1-1)^2 = 0 - \frac{n}{2} + 1 = 1 + \frac{n}{2}$.
But then $1+r-r^2-r^3 = 1+r-1-\frac{n}{2} - r - \frac{n}{2}=0$, and so $\C_2$ is nonempty, a contradiction.
This completes the proof.
\end{proof}


\section{Proof of Theorem~\ref{the:thetheorem}}
\label{sec:proof}

\indent 
We are now ready to finally prove Theorem~\ref{the:thetheorem}.\bigskip

\noindent
{\sc Proof of Theorem~\ref{the:thetheorem}:}
Suppose first that $X$ is a tightly attached half-arc-transitive graph of valency $4$ and
even radius $n$. Theorem~\ref{the:Xgraphs} implies that $X \iso \X_e(m,n;r,t)$, where $m \geq 4$ is even,
$r \in \ZZ_n^*$, $t \in \ZZ_n$ are such that 
$r^m = 1$, $t(r-1) = 0$ and $1 + r + \cdots + r^{m-1} + 2t = 0$. Now Proposition~\ref{pro:arc-tr}
implies that $r^2 \neq \pm 1$ and Lemma~\ref{le:m=6} implies that the condition $(ii)$ of 
the theorem also cannot hold.

Suppose now that $X \iso \X_e(m,n;r,t)$, where $m \geq 4$ is even,
$r \in \ZZ_n^*$, $t \in \ZZ_n$ are such that 
$r^m = 1$, $t(r-1) = 0$ and $1 + r + \cdots + r^{m-1} + 2t = 0$, and none of the two 
conditions of the theorem is fulfilled. Depending on $m$ apply one of
Lemmas~\ref{le:general},~\ref{le:m=8},~\ref{le:m=6} and~\ref{le:m=4}. \hfill $\Qed$\bigskip

The next proposition determines which pairs of the half-arc-transitive graphs $\X_e(m,n;r,t)$
are isomorphic. It transpires that the only possible isomorphisms are those given by
Proposition~\ref{pro:theisomorphisms}. 

\begin{proposition}
\label{pro:alliso}
Let $X = \X_e(m,n;r,t)$, where $m,n \geq 4$ are even integers and $r\in \ZZ_n^*$, 
$t \in \ZZ_n$ are such that  (\ref{eq:Xparamcond}) holds. If $X$ is half-arc-transitive,
then $X \iso \X_e(m',n';r',t')$ if and only if 
$m' = m$, $n' = n$ and one of the following holds:
\begin{itemize}
\item[(i)] $r' = r$ and $t' = t$;
\item[(ii)] $r' = -r$ and $t' = t + r + r^3 + \cdots + r^{m-1}$;
\item[(iii)] $r' = r^{-1}$ and $t' = t$;
\item[(iv)] $r' = -r^{-1}$ and $t' = t + r + r^3 + \cdots + r^{m-1}$.
\end{itemize}
\end{proposition}

\begin{proof}
Proposition~\ref{pro:theisomorphisms} implies that each of the four conditions is sufficient for the
isomorphism to exist. Suppose now that $X \iso \X_e(m',n';r',t') = X'$ and let
$\varphi : X \to X'$ denote one of the isomorphisms. 
Therefore $X$ and $X'$ are both
half-arc-transitive. Fix orientations $D$ and $D'$ of the edges of $X$ and $X'$, respectively, implied
by the half-arc-transitive action of their automorphism groups, and denote the vertex sets of $X$ and $X'$ by 
$\{u_i^j\ |\ i \in \ZZ_m,\ j \in \ZZ_n\}$ and $\{v_i^j\ |\ i \in \ZZ_{m'},\ j \in \ZZ_{n'}\}$, respectively, 
with edges as usual. 

We claim that $\varphi$ either preserves the orientation of every edge or inverts the orientation of every edge.
Suppose on the contrary that for some $w_1,w_2,w_3,w_4 \in V(X)$ and $w_1', w_2', w_3', w_4' \in V(X')$
such that $(w_1,w_2)$, $(w_3,w_4)$, $(w_1',w_2')$ and $(w_3',w_4')$ are arcs in the respective
oriented graphs and that $w_1\varphi = w_1'$, $w_2\varphi = w_2'$, $w_3\varphi = w_4'$ and $w_4\varphi = w_3'$.
There exists some $\alpha \in \Aut X$ mapping $(w_1,w_2)$ to $(w_3,w_4)$ and 
similarly there exists some $\beta \in \Aut X'$ mapping $(w_1',w_2')$ to $(w_3',w_4')$. But then 
$\varphi\beta\varphi^{-1}\alpha^{-1}$ is an automorphism of $X$ interchanging adjacent vertices $w_1$ and $w_2$,
which contradicts the half-arc-transitivity of $X$. Thus our claim holds.

This implies that $\varphi$ maps alternating cycles of $X$ to alternating cycles of $X'$, so since
the lengths of the respective alternating cycles are $2n$ and $2n'$, we have $n' = n$. Thus also $m' = m$.
We can of course assume that $u_0^0\varphi = v_0^0$ and that either 
$u_1^1\varphi = v_1^1$ or $u_1^1\varphi = v_{m-1}^{-r'^{m-1}-t'}$ (otherwise take $\tau\varphi$, where
$\tau \in \Aut X$ is as in Section~\ref{sec:notation}). 

Suppose first that $u_1^1\varphi = v_1^1$. Then in view of the fact that
the alternating cycle of $X$ containing $(u_0^0, u_1^1)$ is mapped to the 
alternating cycle of $X'$ containing $(v_0^0, v_1^1)$, we clearly have 
$u_0^i\varphi = v_0^i$ and $u_1^i\varphi = v_1^i$ for all
$i \in \ZZ_n$. Therefore, $u_1^r\varphi = v_1^r$ and $u_1^0\varphi = v_1^0$, so since
$u_1^r$ and $u_1^0$ have a unique common neighbor $u_2^r$, so do $v_1^r$ and $v_1^0$. Moreover,
this common neighbor is a vertex of the form $v_2^i$. Thus either $r -r' = 0$ or $r + r' = 0$, that is
$r' = r$ or $r' = -r$.

Suppose now that $u_1^1\varphi = v_{m-1}^{-r'^{m-1}-t'}$. Then $u_0^1\varphi = v_0^{-r'^{m-1}}$, 
and so a similar argument as in the previous paragraph shows that 
$u_1^r\varphi = v_{m-1}^{-rr'^{m-1}-t'}$ and $u_1^0\varphi = v_{m-1}^{-t'}$.
Thus either $-rr'^{m-1}-t'-r'^{m-2} = -t'$ or $-rr'^{m-1}-t'+r'^{m-2} = -t'$. In the first case
$r'^m = 1$ implies $r = -r'^{m-1}$, and so $r' = -r^{-1}$, and in the second case $r' = r^{-1}$.

Using Proposition~\ref{pro:theisomorphisms} it now suffices to show that $t' = t$ whenever $r' = r$.
Note that $r' = r$ implies (in view of $1+r+\cdots + r^{m-1}¸+ 2t = 0$) that $2t = 2t'$, so
either $t' = t$ or $t' = t +\frac{n}{2}$. Moreover, since $r^2 \neq \pm 1$, 
the previous paragraphs of this proof
imply that we have $u_0^i\varphi = v_0^i$ and $u_1^i\varphi = v_1^i$ for all $i \in \ZZ_n$.
The vertices $u_1^i$ and $u_1^{i+r}$ have a unique (recall that $r \neq \pm 1$ in view of half-arc-transitivity)
neighbor $u_2^{i+r}$, and so since $v_2^{i+r}$ is the unique common neighbor of $v_1^i$ and $v_1^{i+r}$,
we have $u_2^{i}\varphi = v_2^i$. Continuing inductively we get $u_i^j \varphi = v_i^j$ for all
$i \in \ZZ_m$, $j \in \ZZ_n$. Since $u_{m-1}^0$ is adjacent to $u_0^t$, the vertex
$v_{m-1}^0$ is adjacent to $v_0^t$, and so either $t = t'$ or $t = r^{m-1}+t'$. But the latter case
is impossible in view of the fact that either $t' = t$ or $t' = t +\frac{n}{2}$ (recall that $r \in \ZZ_n^*$). 
Thus $t' = t$, as claimed.
\end{proof}\bigskip

We end this paper by the following observations. In \cite{MS98} \v Sajna considered metacirculants
$M(r;m,n)$, which are defined to have vertex set 
$\{u_i^j\ |\ i \in \ZZ_m,\ j \in \ZZ_n \}$, where $r \in \ZZ_n^*$ is such that $r^m = \pm 1$,
and edge set $\{u_i^j u_{i+1}^{j \pm r^i}\ |\ i \in \ZZ_m,\ j \in \ZZ_n\}$. It turns out that
the graphs $M(r;m,n)$ admit a half-arc-transitive group action. One of \v Sajna's
goals was to determine which of them are half-arc-transitive and which are arc-transitive.
The case that remained unsolved was the case 
$n \equiv 0 \pmod 4$, and $r$ of order $4$. In this case $M(r;4,n)$
consists of two isomorphic connected components,
one containing vertices $u_i^j$, for which $i \equiv j \pmod 2$, and the other containing vertices 
$u_i^j$, for which $i \not\equiv j \pmod 2$.
The next corollary of Lemma~\ref{le:m=4} determines which of the above metacirculants are half-arc-transitive
and which are arc-transitive, thus completing the work initiated by \v Sajna.

\begin{corollary}
\label{cor:sajna}
Let $n_1$ be even, let $n = 2n_1$ and let $r \in \ZZ_n^*$ be of order $4$. Let $X$ be one of the two
connected components of $M(r;4,n)$. 
Then $X$ is half-arc-transitive unless $r^2 \equiv \pm 1 \pmod {n_1}$ in which case
it is arc-transitive.
\end{corollary}

\begin{proof}
In view of \cite[Lemma~3.7.]{MS98} we have $M(r;4,n) \iso M(-r;4,n)$, and so we can assume $r < n_1$. We can thus 
regard $r$ also as an element of $\ZZ_{n_1}$. Denote the
corresponding element of $\ZZ_{n_1}$ by $r'$.
Let $t \in \ZZ_n$, $t < n_1$, be such that $-1-r-r^2-r^3 = 2t$. Since $r \in \ZZ_n^*$ and
$n$ is even, such $t$ exists and is unique. Denote the corresponding element of $\ZZ_{n_1}$
by $t'$.
Since the two components of $M(r;4,n)$ are isomorphic we can assume that $X$ contains $u_0^0$.
We show that $X \iso \X_e(4,n_1;r',t')$ and that $r', t'$ satisfy
conditions $(\ref{eq:Xparamcond})$. Then Lemma~\ref{le:m=4} implies that
$X$ is half-arc transitive unless $r'^2 = \pm 1$, in which case it is arc-transitive by  
Proposition~\ref{pro:arc-tr}.

In the following paragraph the elements of $\ZZ_{n_1}$ will constantly be viewed as the respective
elements of $\ZZ_n$. This should cause no confusion. 
Denote the vertices of $\X_e(4,n_1;r',t')$ by $\{v_i^j\ |\ i \in \ZZ_4,\ j \in \ZZ_{n_1}\}$, with edges as usual.
The isomorphism $\varphi : X \to \X_e(4,n_1;r',t')$ is defined as follows. 
For each vertex of the form $u_0^j$ of $X$ (note that $j$ is even) there is a unique $j_1 \in \ZZ_{n_1}$, such that
$2j_1 = j$ in $\ZZ_n$. Set $u_0^j\varphi = v_0^{j_1}$. For each vertex of the form $u_1^j$ of $X$ (note that
$j$ is now odd) there is a unique
$j_1 \in \ZZ_{n_1}$, such that $j = 2j_1 - 1$ in $\ZZ_n$. Set $u_1^j\varphi = v_1^{j_1}$.
For each $u_2^j$ there is a unique $j_1 \in \ZZ_{n_1}$, such that $j = 2j_1 - 1 - r$ in $\ZZ_n$. 
Set $u_2^j\varphi = v_2^{j_1}$. Finally, for each $u_3^j$ there is a unique $j_1 \in \ZZ_{n_1}$, such
that $j = 2j_1 - 1 - r - r^2$ in $\ZZ_n$. 
Set $u_3^j\varphi = v_3^{j_1}$. Clearly $\varphi$ is a bijection. 
It is easy to see that it also preserves the adjacencies. We leave this to the reader.
Clearly $r'^4 = 1$. Moreover, $1+r'+r'^2+r'^3+2t' = 0$ in $\ZZ_{n_1}$ since $1+r+r^2+r^3+2t = 0$ in $\ZZ_n$.
Note that the latter equation and $r^4 = 1$ in $\ZZ_n$ also imply 
$r+r^2+r^3+1+2tr = 0$, and so $2t(r-1) = 0$ in $\ZZ_n$. 
But then $t'(r'-1) = 0$ in $\ZZ_{n_1}$, which completes the proof.
\end{proof}


\end{document}